\documentclass{amsart}

\usepackage{amsmath,amsthm,amssymb,latexsym,mathtools}

\usepackage[T1]{fontenc}
\RequirePackage{fix-cm}

\usepackage[dvipsnames]{xcolor}
\usepackage{graphicx}
\usepackage{tikz}
\usetikzlibrary{decorations.pathreplacing}
\usepackage{float}
\usepackage{standalone}

\usepackage{xspace}
\usepackage[all]{xy}
\usepackage{pinlabel}
\usepackage{enumerate}
\usepackage[font=small,labelfont=bf]{caption}
\usepackage{subcaption}
\usepackage{wrapfig}
\usepackage[percent]{overpic}
\usepackage{stackrel}
\usepackage{accents}
\usepackage{multirow}
\usepackage{circuitikz}
\usepackage{lipsum}
\usepackage{comment}

\usepackage{hyperref}
\hypersetup{
    colorlinks=true,
    linkcolor=MidnightBlue,
    citecolor=OliveGreen,
    urlcolor=BrickRed,
    pdfborder={0 0 0} 
}
\usepackage[capitalize]{cleveref}

\usepackage{tikz}
\usetikzlibrary{decorations.pathreplacing}
\usepackage{float}
\usepackage{standalone}
 \usepackage{graphicx}
  \graphicspath{{pictures/}}
  \usepackage{xspace}
  \usepackage[all]{xy}
  \usepackage{pinlabel}
  \usepackage{enumerate}
  \usepackage[font=small,labelfont=bf]{caption}
  \usepackage{wrapfig}
\usepackage{lipsum}
\usepackage[percent]{overpic}
\usepackage{caption}
\usepackage{subcaption}
\usepackage{stackrel}
\usepackage{accents}
\usepackage{multirow}

  \renewcommand{\AA}{\mathbb{A}}
  
  \newcommand{\CC}{\mathbb{C}}
  \newcommand{\DD}{\mathbb{D}}

  \newcommand{\HH}{\mathbb{H}}

  \newcommand{\QQ}{\mathbb{Q}}
  \newcommand{\RR}{\mathbb{R}}

  \newcommand{\figref}[1]{Figure~\ref{#1}}

\theoremstyle{definition}
\newtheorem{proposition}{Proposition}[section]
\newtheorem{corollary}[proposition]{Corollary}
\newtheorem{lemma}[proposition]{Lemma}
\newtheorem{theorem}[proposition]{Theorem}
\newtheorem{definition}[proposition]{Definition}
\newtheorem{notation}[proposition]{Notation}

\newtheorem{remark}[proposition]{Remark}

\crefname{theorem}{theorem}{theorems}
\Crefname{theorem}{Theorem}{Theorems}
\crefname{proposition}{proposition}{propositions}
\Crefname{proposition}{Proposition}{Propositions}
\crefname{lemma}{lemma}{lemmas}
\Crefname{lemma}{Lemma}{Lemmas}
\crefname{corollary}{corollary}{corollaries}
\Crefname{corollary}{Corollary}{Corollaries}
\crefname{definition}{definition}{definitions}
\Crefname{definition}{Definition}{Definitions}
\crefname{remark}{remark}{remarks}
\Crefname{remark}{Remark}{Remarks}
\crefname{example}{example}{examples}
\Crefname{example}{Example}{Examples}
\crefname{conjecture}{conjecture}{conjectures}
\Crefname{conjecture}{Conjecture}{Conjectures}

\title[Moduli spaces of open strings have polylog volumes]
{Moduli spaces of open strings have polylogarithmic Mirzakhani volumes}
\author{Yi Huang}

\address{Yi Huang, Room C651 Shuangqing Complex Building A, Yau Mathematical Sciences Center,
Tsinghua University, Haidian District
Beijing 100084, China}
 \email{yihuangmath@tsinghua.edu.cn}

\author{Ivan Telpukhovskiy}

\address{Ivan Telpukhovskiy, Room C656 Shuangqing Complex Building A, Yau Mathematical Sciences Center,
Tsinghua University, Haidian District
Beijing 100084, China}
 \email{ivantelp@tsinghua.edu.cn}

\date{\today}

\begin{document}
\begin{abstract}
We show that the Mirzakhani volume, as introduced by Chekhov, of the moduli space of every crowned hyperbolic surface is naturally expressible as a sum of Gaussian rational multiples of polylogarithms evaluated at $\pm1$ and $\pm\sqrt{-1}$.
\end{abstract}
\maketitle

\section{Introduction}

In \cite{mirz_simp}, Mirzakhani showed that the Weil--Petersson volume of the moduli space $\mathcal{M}_{g,m}(b_1,\ldots,b_m)$ of hyperbolic structures on a genus $g$ surface with $m$  boundaries, with specified boundary lengths $b_1,\ldots,b_m$, is a polynomial in $\mathbb{Q}_{\geqslant 0}[\pi^2,b_1^2,\ldots,b_m^2]$. She subsequently reinterpreted the coefficients of these polynomials via Duistermaat--Heckman theory as intersection numbers on moduli spaces of pointed curves, and used this to obtain a novel proof of Witten's conjecture \cite{witten_conjecture} --- also known as Kontsevich's theorem \cite{MR1171758}, that the generating function for the aforementioned intersection numbers should satisfy the Korteweg-de Vries hierarchy.\medskip

There has been great progress in building an analogue of Witten's conjecture, at least from the algebraic geometric perspective \cite{zbMATH06877694,  zbMATH06824078, zbMATH07023773, zbMATH06576752, zbMATH06497098, buryak2017matrix, zbMATH07951658}, for intersection numbers of moduli spaces of \emph{open strings} --- i.e.: Riemann surfaces with marked points (or equivalently, punctures) on their boundaries. However, any attempt to find a parallel for Mirzakhani's hyperbolic geometric approach to Witten's conjecture for open strings has to contend with the non-finiteness of the Weil--Petersson volumes of the moduli spaces at hand.\medskip

In \cite{chekhov2024}, Chekhov introduces an action on moduli spaces of crowned hyperbolic surfaces (i.e.: uniformizations of open strings) of fixed neck holonomy, and we refer to the integrals of the induced measures over moduli spaces as \emph{Mirzakhani volumes}. Chekhov computes a handful of examples, and we extend Chekhov's work to full generality (\Cref{thm:fixedneckvol-intro}, \Cref{thm:ngonvol-intro}, \Cref{cor:fixedneckvolumegeneral-intro}), and beyond (\Cref{thm:compositethm}, \Cref{rmk:polylogs}).

\subsection{Main results}

\begin{notation}
    We maintain the following notation throughout this paper:
    \begin{itemize}
        \item $\Sigma$ is a topological surface that admits a crowned hyperbolic surface structure (\Cref{notn:sigmasurface});
        \item $X$ denotes a crowned hyperbolic surface (\Cref{defn:crowned})
        \item $\mathcal{M}_\Sigma(\vec{b})$ is the set of crowned surfaces with cuffs of length $\vec{b}\in[0,\infty)^m$ (\Cref{def:fixed-holonomy}), and $\mathcal{M}_\Sigma(\vec{b}|\vec{d})\subset \mathcal{M}_\Sigma(\vec{b})$ is the subset of $\mathcal{M}_\Sigma(\vec{b})$ consisting of crowned hyperbolic surfaces with necks of length $\vec{d}\in[0,\infty)^l$ (\Cref{defn:neck} and \Cref{def:fixed-holonomy});
        \item $V_\Sigma(\vec{b})$ and $V_\Sigma(\vec{b}|\vec{d})$ respectively denote the Mirzakhani volumes (\Cref{defn:mirzvolume}) of $\mathcal{M}_\Sigma(\vec{b})$ and $\mathcal{M}_\Sigma(\vec{b}|\vec{d})$.
    \end{itemize}
\end{notation}

\begin{notation}[disks and annuli]\label{notn:specialsurfaces}
    We especially write
    \begin{itemize}
    \item $\Sigma=\mathbb{D}_n$ for an \emph{$n$-gon}, i.e.: a disk with $n\geqslant 3$ boundary punctures.  
    \item $\Sigma=\mathbb{A}_{n}$ for an \emph{$n$-crown}, i.e.: a closed annulus with $n\geqslant 1$ boundary punctures on one boundary and no punctures on the other boundary.
    \item $\Sigma=\mathbb{A}_{a_1,a_2}$ for an \emph{$(a_1,a_2)$-annulus}, i.e.: a closed annulus with $a_1$ punctures on one boundary and $a_2$ boundary punctures on other boundary.    
    \end{itemize}
    To clarify: we regard the cuff/neck of $\mathbb{A}_n$ (and also $\mathbb{A}_{a_1,a_2}$) as a neck rather than a cuff.
\end{notation}

Chekhov computes the Mirzakhani volume of $\mathcal{M}_{\mathbb{A}_{3}}(d)$ in \cite[2.5.1]{chekhov2024}. We generalize this for moduli spaces of all $n$-crowns with fixed neck lengths, which in turn yields the Mirzakhani volumes for moduli spaces of all crowned hyperbolic surfaces (with fixed neck lengths) apart from moduli spaces of $n$-gons:

\begin{theorem}
[\Cref{thm:fixedneckvol}]
\label{thm:fixedneckvol-intro}
The Mirzakhani volume $V_{\mathbb{A}_n}(d)$ of $\mathcal{M}_{\mathbb{A}_n}(d)$ is:
\begin{align*}
V_{\mathbb{A}_n}(d) 
=
\left\{
\begin{array}{rl}
\frac{d}{\sinh{d/2}} \cdot 
\frac{\prod_{j=1}^{k-1} (d^2+(2j)^2\pi^2)}{2(n-1)!},
&\quad\text{if }n=2k\text{ is even;}
\\[0.25em]
\frac{1}{\cosh{d/2}} \cdot  \frac{\prod_{j=1}^{k} (d^2+(2j-1)^2\pi^2 )}{2(n-1)!},
&\quad\text{if }n = 2k+1\text{ is odd.}
\end{array}
\right.
\end{align*}

\end{theorem}

\begin{corollary}[\Cref{cor:fixedneckvolumegeneral}]
\label{cor:fixedneckvolumegeneral-intro}
    For $\Sigma=\mathbb{A}_{a_1,a_2}$, the Mirzakhani volume of $\mathcal{M}_{\mathbb{A}_{a_1,a_2}}(d)$ is
    \[
    V_{\mathbb{A}_{a_1,a_2}}(d)
    =
    d\cdot V_{\mathbb{A}_{a_1}}(d)\cdot V_{\mathbb{A}_{a_2}}(d).
    \]
    For $\Sigma\neq\mathbb{D}_n,\mathbb{A}_n,\mathbb{A}_{a_1,a_2}$, let $g$ denote the genus of $\Sigma$. The volume of $\mathcal{M}_\Sigma(\vec{b}|\vec{d})$ is
    \[
    V_\Sigma(\vec{b}|\vec{d})
    =
    V_{g,m+l}(\vec{b},\vec{d})
    \cdot
    \prod_{k=1}^{l}d_k \cdot V_{\mathbb{A}_{a_k}}(d_k),
    \]
    where $V_{g,m+l}(\vec{b},\vec{d})\in\mathbb{Q}_{>0}[\pi^2,b_1^2,\ldots,b_m^2,d_1^2,\ldots,d_l^2]$ is the Weil--Petersson volume of the moduli space $\mathcal{M}_{g,m+l}(\vec{b},\vec{d})$.
\end{corollary}

We further determine the Mirzakhani volume of the moduli space $\mathcal{M}_\Sigma(\vec{b})$, where there are no neck length constraints. We obtain the following composite theorem as the highlight of this paper:

\begin{theorem}[\Cref{thm:ngonvol}, \Cref{thm:crownfullvol}, \Cref{thm:vol-biannuli}, \Cref{thm:generalmirzvolumes}]
\label{thm:compositethm}
For any $\Sigma$, the volume of the $\mathcal{M}_\Sigma(\vec{b})$, where $\vec{b}=\varnothing$ when there are no cuffs, satisfies
\begin{align*}
V_\Sigma(\vec{b}) \in \QQ_{>0}[b_1^2,\ldots,b_m^2, \log2, \zeta(j),\beta(2k)],
\quad
j,2k \in \{2,3,\ldots,\dim_{\mathbb{R}}\mathcal{M}_{\Sigma}(\vec{b})\},
\end{align*}
where $\zeta(s)$ is the Riemann zeta function and $\beta(s)$ is the Dirichlet beta function. Moreover, the volume is a homogeneous polynomial with rational coefficients, where $b_i$ is degree~1, $\log 2$ is degree $1$, $\zeta(j)$ is degree $j$, and $\beta(2k)$ is degree $2k$.
\end{theorem}

\begin{remark}
\label{rmk:polylogs}
The Mirzakhani volumes  for $\mathcal{M}_\Sigma(\vec{b})$ are naturally expressible as polylogarithms of $\pm1$ and $\pm\sqrt{-1}$:
\begin{align*}
\log 2=-\text{Li}_1(-1),\quad
\zeta(i)=\text{Li}_i(1),\quad
\beta(2j)=\tfrac{\sqrt{-1}}{2}
\left(\text{Li}_{2j}(-\sqrt{-1})-\text{Li}_{2j}(\sqrt{-1})\right).
\end{align*}

\end{remark}

\subsubsection{Mirzakhani volumes of moduli spaces of $n$-gons}

Chekhov determined the $n=5,6$ cases in \cite[6.2.1 and 6.2.2]{chekhov2024}, and in the first version of this paper, we computed the $n=7,8$ cases and found a conjectural formula based on numerical calculations and established an upper bound for the volume for general $n$. Also in response to this earlier version, Timothy Budd and Leonid Chekhov, who had independently obtained progress on these volume calculations, kindly shared with us that the desired integral has been derived on MathOverflow \cite{MO402985}, which we have been able to use to prove:

\begin{theorem}[\Cref{thm:ngonvol}]
\label{thm:ngonvol-intro}
The Mirzakhani volume of the moduli space of ideal $n$-gons, for $n\geqslant 3$, is
\begin{align*}
V_{\DD_n} = \left\{
\begin{array}{rl}
\frac{2}{n-2}\cdot\frac{(n-4)!!}{(n-3)!!}\pi^{n-4},
&\quad\text{if }n\text{ is even,}
\\[0.25em]
  \frac{1}{n-2}\cdot\frac{(n-4)!!}{(n-3)!!}\pi^{n-3},
&\quad\text{if }n\text{ is odd.}
\end{array}
\right.    
\end{align*}    
\end{theorem}

\begin{remark}
This theorem is equivalent (\Cref{prop:genfunction}) to the statement that the generating function for Mirzakhani volumes of moduli spaces of $n$-gons is  
\begin{align*}
\frac{\arcsin{\pi x}}{\pi x}
+
x\left(\frac{\arcsin{\pi x}}{\pi x}\right)^2.
\end{align*}
\end{remark}

\begin{remark}[potential connections to Ap\'ery constants]
    There appears to be a curious agreement between the highest index non-zero Ap\'{e}ry numbers for Grassmannians $\mathrm{Gr}(2,n)$ and $V_{\mathbb{D}_n}$ (see the table on page~3 of Galkin's paper \cite{galkin_apery}). This was first communicated to us by Zichang Wang\footnote{Confusingly, there's another excellent student of the same transliterated name in the same year at the same college who is working on mathematics with physical motivations, and whose Chinese names are only one character apart.} --- a (then-undergrad) student at Qiuzhen College, Tsinghua University.
\end{remark}

\subsection*{Acknowledgements} The authors wish to thank Enoch Guo for assistance with some of the computer-based volume computations for this paper, and Timothy Budd and Leonid Chekhov for pointing us to \cite{MO402985}. The second listed author gratefully acknowledges support from the Beijing Natural Science Foundation (International Scientists Project), Funding No. 1S24065.

\section{Background}

We introduce background and notation for defining Chekhov's notion of Mirzakhani volumes.

\subsection{Open strings and crowned hyperbolic surfaces}

From the complex analytic perspective, an \emph{open string} or an \emph{open Riemann surface} with punctures is a Riemann surface with (smooth) boundary with interior and boundary punctures \cite[Section~1.3]{zbMATH07951658}. Standard Schottky doubling arguments (see, e.g.: \cite[Pg.~24]{huang_thesis}) show that such surfaces uniformize to a complete finite-area hyperbolic surface with bi-infinite geodesic boundary components --- i.e.: a crowned hyperbolic surface where every end is either a cusp or a tine (half-cusp). We generalize very mildly beyond this context to align with what Mirzakhani did: we allow also for non-cuspidal cuffs. In terms of complex geometry, these are Riemann surfaces minus finitely many interior points, boundary points, and (closed) disks in the interior of the surface. 

\begin{notation}
\label{notn:sigmasurface}
For the remainder of the article let $\Sigma:=\Sigma_{m,\vec{a}}$, with $m\in\mathbb{Z}_{\geqslant 0}$ and $\vec{a}=(a_1,\ldots,a_l)\in\mathbb{Z}_{>0}^l$, denote a bordered topological surface of finite type, where
\begin{itemize}
\item
$n=\sum_{k=1}^{l}a_k$ boundary punctures, which we label $\{q_1,\ldots,q_n\}$; 
\item
$\Sigma\cup\{q_1,\ldots,q_n\}$ is a surface with $m$ interior punctures, which we label $\{p_1,\ldots,p_m\}$, and $l\geq0$ boundary curves, which we label $\{\delta_1,\ldots,\delta_l\}$; 
\item
the boundary curves $\delta_i$ respectively contain $a_i\geqslant 1$ punctures. We order the punctures $\{q_1,\ldots,q_n\}$ so that $q_1,\ldots,q_{a_1}\in \delta_1$, and $q_{a_1+1},\ldots q_{a_1+a_2}\in \delta_2$, and so forth, and so that successively labelled punctures appear consecutively along each boundary. 
\end{itemize}
We further impose the condition that $\Sigma$ admits a hyperbolic metric with geodesic boundary, meaning that
\begin{itemize}
    \item if $g=0$ and $l=0$, then $m\geqslant3$; 
    \item if $g=0$, $l=1$ and $m=0$, then $n\geqslant3$;
    \item if $g=1$ and $l=0$, then $m\geqslant1$.
\end{itemize}
\end{notation}

\begin{definition}[crowned hyperbolic surface]
\label{defn:crowned}
We refer to $\Sigma$, equipped with a complete hyperbolic metric $h$, as a \emph{crowned hyperbolic surface} if
\begin{itemize}
    \item $(\Sigma,h)$ has bi-infinite geodesic boundaries;
    \item each of the boundary punctures $q_j$ uniformizes to half a cusp.
\end{itemize}
We henceforth use $X=(\Sigma,h)$ to denote crowned hyperbolic surfaces.
\end{definition}

\subsection{The anatomy of a crowned hyperbolic surface}

\begin{definition}[arches, cuffs, and tines]
\label{defn:anatomy}
Given a crowned hyperbolic surface $X=(\Sigma,h)$, we refer to the bi-infinite geodesic boundaries of $X$ as \emph{arches}, and denote them by $\{\alpha_j\}_{j=1,\ldots,n}$. The following types of behaviour can hold at the punctures of $\Sigma$:
\begin{itemize}
    \item Interior punctures $p_i$ geometrise as flares/funnels or cusps, (any of) which we will refer to as the $i$-th \emph{cuff}. We write $\beta_i$ for a simple closed loop which goes around $p_i$ (once). Whenever possible, we assume that the cuff is the geodesic representative of $\beta_i$.
    \item Boundary punctures $q_j$ geometrise to a ``half-cusp'' called a \emph{tine} (see \cite{huang_thesis}). We refer to $q_j$ as the $j$-th tine.
\end{itemize}
\end{definition}

\begin{definition}[neck]
\label{defn:neck}
    For each boundary curve $\delta_k$ on $X\cup\{q_1,\ldots,q_n\}$, choose a simple isotopy representative $\nu_k$ contained in the interior of $X$. We refer to $\nu_k$ as the $k$-th \emph{neck} of $X$.
\end{definition}

\begin{definition}[cuff and neck length]
    We refer to the hyperbolic geometric lengths of the $i$-th cuff (resp. the $k$-th necks) as the $i$-th cuff length (resp. the $k$-th neck length), where cuffs around cusps have cuff length $0$.  
\end{definition}

\begin{center}
    \begin{figure}
        \centering
        \includegraphics[width=0.5\linewidth]{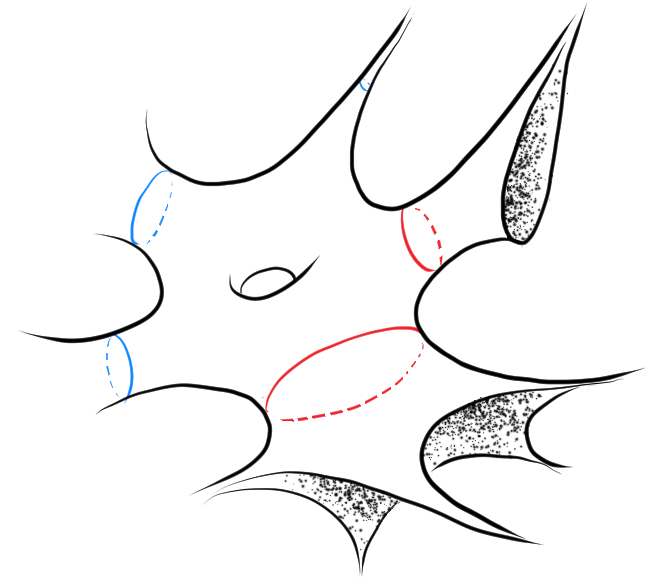}
        \caption{An example of a crowned hyperbolic surface with genus $1$, three {\color{blue}cuffs} and two {\color{red}necks} respectively isolating an $1$-crown and a $5$-crown.}
        \label{fig:crownedsurface}
    \end{figure}
\end{center}

\subsection{Teichm\"{u}ller and moduli spaces of crowned hyperbolic surfaces}

\begin{definition}[Teichm\"uller space]
The \emph{Teichm\"uller space} $\mathcal{T}_\Sigma$ is the space of isotopy classes of crowned hyperbolic surface metrics on $\Sigma$.    
\end{definition}

\begin{definition}[moduli space]
\label{defn:modulispace}
The \emph{moduli space} $\mathcal{M}_\Sigma$ is the space of orientation-preserving isometry classes of crowned hyperbolic surfaces.  
\end{definition}

The (pure) mapping class group $\mathrm{MCG}_\Sigma$, consisting of isotopy classes of puncture-preserving homeomorphisms of $\Sigma$, is the quotient group relating $\mathcal{T}_\Sigma$ and $\mathcal{M}_\Sigma$:
\[
\mathcal{M}_\Sigma=\mathcal{T}_\Sigma/\mathrm{MCG}_\Sigma.
\]

\begin{definition}[moduli spaces with fixed holonomy]
\label{def:fixed-holonomy}
For $\vec{b}=(b_1,\ldots,b_m)\in[0,\infty)^m$, define
\begin{align*}
\mathcal{T}_\Sigma(b_1,\ldots,b_m)
&=:\mathcal{T}_\Sigma(\vec{b})
\subset\mathcal{T}_\Sigma
\quad\text{and}\\
\mathcal{M}_\Sigma(b_1,\ldots,b_m)
&=:\mathcal{M}_\Sigma(\vec{b})
\subset\mathcal{M}_\Sigma,
\end{align*}
to respectively be the subsets of crowned hyperbolic surfaces in $\mathcal{T}_\Sigma$ and $\mathcal{M}_\Sigma$ whose $i$-th cuff $\beta_i$ is of length $b_i$ for all $i=1,\ldots,m$. Further, given $\vec{d}=(d_1,\ldots,d_l)\in[0,\infty)^l$, define
\begin{align*}
\mathcal{T}_\Sigma(b_1,\ldots,b_m|d_1,\ldots,d_l)
&=:\mathcal{T}_\Sigma(\vec{b}|\vec{d})
\subset\mathcal{T}_\Sigma(\vec{b})
\quad\text{and}\\
\mathcal{M}_\Sigma(b_1,\ldots,b_m|d_1,\ldots,d_l)
&=:\mathcal{M}_\Sigma(\vec{b}|\vec{d})
\subset\mathcal{M}_\Sigma(\vec{b}).
\end{align*}
to respectively be the subsets of hyperbolic surfaces whose $k$-th neck lengths is $d_k$ for all $k=1,\ldots,l$.
\end{definition}

\begin{notation}[fixed holonomy for special surfaces]
When $\Sigma=\mathbb{D}_n$, $\mathbb{A}_{n}$, or $\mathbb{A}_{a_1,a_2}$, the above notation for Teichm\"uller and moduli spaces with fixed holonomy may not make sense due to the absence of cuffs and/or necks. Instead, we write:
    \begin{itemize}
        \item  $\mathcal{T}_{\mathbb{D}_n}$ and $\mathcal{M}_{\mathbb{D}_n}$;
        \item $\mathcal{T}_{\mathbb{A}_{n}}(d)$ and $\mathcal{M}_{\mathbb{A}_{n}}(d)$;
        \item $\mathcal{T}_{\mathbb{A}_{a_1,a_2}}(d)$ and $\mathcal{M}_{\mathbb{A}_{a_1,a_2}}(d)$.
    \end{itemize}
\end{notation}

\begin{remark}
    Teichm\"uller spaces with fixed cuff lengths are open balls (possibly of dimension $0$). In particular, $\mathcal{T}_\Sigma(\vec{b})$ embeds as a half-dimensional subspace in 
    \[
    \mathcal{T}_{D\Sigma}(b_1,\ldots,b_m,b_1,\ldots,b_m,\underbrace{0,\ldots,0}_{n\text{ terms }}),
    \]
    and hence has dimension  $\frac{1}{2}(3|\chi(D\Sigma)|-2m-n)$. The fixed (cuff and) neck length subset $\mathcal{T}_\Sigma(\vec{b}|\vec{d})$ embeds as an $l$ co-dimensional subset of $\mathcal{T}(\vec{b})$ and hence has dimension $ \frac{1}{2}(3|\chi(D\Sigma)|-2m-2l-n)$.
\end{remark}

\subsection{Shearing coordinates on Teichm\"uller space of crowned surfaces}
\label{subsec:shears}

The contents of this subsection is well-known to experts. See, for example, \cite[Section~7.4]{martelli2022introductiongeometrictopology} for a definition of the shearing parameter between two adjacent ideal hyperbolic triangles and the construction of shearing coordinates on Teichm\"uller space of surfaces with punctures.

\begin{notation}[shearing sign convention]
    We take \emph{left} shears to be positive.
\end{notation}

\begin{definition}[ideal triangulation]
Given a surface $\Sigma=\Sigma_{m,\vec{a}}$, we call a collection $\triangle$ of disjoint simple (open) arcs in $\Sigma$ with (only) endpoints placed at the interior and boundary punctures of $\Sigma$ an \emph{ideal triangulation} of $\Sigma$ if
\begin{itemize}
    \item no arc in $\triangle$ bounds an open disk, 
    \item no arc is peripheral in the sense of being isotopic to one of the arches,
    \item no two arcs in $\triangle$ are isotopic via isotopy fixing the endpoints, and
    \item $\triangle$ is maximal in cardinality among sets satisfying above conditions.
\end{itemize}
We identify ideal triangulations up to isotopies of $\Sigma$.  
\end{definition}

\begin{lemma}
Given a surface $\Sigma=\Sigma_{m,\vec{a}}$, 
with $m$ interior punctures and $l$ boundary curves with $n=\sum_{k=1}^l a_k$ boundary punctures so that there are $a_k$ punctures on the $k$-th boundary (see \Cref{notn:sigmasurface}), every ideal triangulation $\triangle$ of $\Sigma$ has cardinality
\[
|\triangle|=\tfrac{1}{2}(3|\chi(D\Sigma)|-n),
\]  
where $D\Sigma$ is the surface obtained by doubling $\Sigma$ along its boundary. Furthermore, if the genus of $\Sigma$ is $g$, then
\[
|\triangle|=6g-6+3m+3l+n.
\]
\end{lemma}
\begin{proof}
Let $\Sigma'$ denote the opposite-oriented copy of $\Sigma$ in $D\Sigma$. Then the ideal triangulation $\triangle$ and its reflection $(\triangle')$ in $\Sigma'$, combined with the boundary arches $\{\alpha_1,\ldots,\alpha_n\}$, give an ideal triangulation of $D\Sigma$. Any ideal triangulation of $D\Sigma$ has cardinality $3|\chi(D\Sigma)|$, and hence $
|\triangle|=\frac{1}{2}(3|\chi(D\Sigma)|-n)
$. Now, if the genus of $\Sigma$ is $g$, then the genus of $D\Sigma$ is $2g+l-1$, while the number of (interior) punctures of $D\Sigma$ is $n+2m$, and hence $\chi(D\Sigma) = 2-2(2g+l-1)-(n+2m) = 4 - 4g -2m -2l-n$. It follows that $\frac{1}{2}(3|\chi(D\Sigma)|-n) = 6g-6+3m+3l+n$.
\end{proof}

Given a crowned hyperbolic surface $X \in \mathcal{T}_{\Sigma}(\vec{b})$, there is a unique geodesic representative of $\triangle=\triangle(X)$ such that if $b_i>0$ (i.e. $p_i$ is geometrised as a flare), then the arcs in $\triangle$ with $p_i$ as an endpoint geometrize to geodesics that spiral towards the $i$-th cuff in the following way: from the perspective of a point on the cuff and looking away from the funnel, the spiralling geodesics pass to the left.

\begin{proposition}[shearing coordinates, see e.g. {\cite[Cor~2.22]{huang_thesis}}]
\label{prop:shearingcoords}
Let  $\triangle$ be an ideal triangulation of $\Sigma$, and let $s_i$ denote the shearing parameter associated to the $i$-th arc of $\triangle$. Then the map 
\[
s_{\triangle} : \mathcal{T}_{\Sigma}(\vec{b}) \to\RR^{\triangle} 
\]
is a real-analytic embedding into the affine subspace in $\RR^{\triangle} $ given by the following conditions for $i=1,\dotsc,m$: the collection of $\{s_j\}$ corresponding to arcs incident to $\beta_i$ sum (with multiplicity) to $b_i$.
\end{proposition}

\begin{remark}
    Strictly speaking, the fact that the $\{s_j\}$ are not independent means that we should only include a subset of these shearing parameters to get a true coordinate chart. We nevertheless refer to $s_\triangle$ as shearing coordinates.
\end{remark}

\subsection{Lambda length coordinates}

We will refer to one more coordinate system related to Teichm\"uller spaces called lambda length coordinates. These coordinates parametrize the (horocycle-)decorated Teichm\"uller space \cite{pennercoords}, as opposed to the Teichm\"uller space itself, and each $\lambda$-length is geometric invariant assigned to a geodesic in $\HH$ that is decorated by (i.e.: paired with) horocycles at each of its ends.

\begin{definition}[lambda length]
Given a geodesic $\gamma\subset\HH$ and a pair of horocycles $h_1, h_2$ whose (distinct) centers are the endpoints of $\gamma$, let $d$ denote the signed hyperbolic distance along $\gamma$ between the points $h_1 \cap \gamma$ and $h_2 \cap \gamma$, where the sign of $d$ is taken to be positive if and only if $h_1$ and $h_2$ are disjoint. The \emph{$\lambda$-length} or \emph{lambda length} of $h_1, h_2$ is
\begin{align}
\label{eq: def lambda length}
\lambda(h_1,h_2) = e^{d/2}.
\end{align}
\end{definition}

\begin{definition}[cusp/tine decoration]
We say a cusp $\beta_i$ (a tine $q_j$) of a crowned hyperbolic surface $X$ is \textit{decorated} if a horocycle centered at $\widetilde{\beta}_i \in \partial_\infty \widetilde{X} \subset \partial \HH \,\,(\widetilde{q}_j \in \partial_\infty \widetilde{X} \subset \partial \HH)$ is chosen in a $\pi_1(X)$-equivariant way. A crowned surface is called a \emph{decorated crowned surface} if all its cusps and tines are decorated; we denote it by $(X,\textbf{h})$, where $\textbf{h}$ denotes the collection of decorating horocycles.
\end{definition}

\begin{definition}[decorated Teichm\"uller space and moduli space]
    The \emph{decorated Teichm\"uller space} $\widetilde{\mathcal{T}}_{\Sigma}$ is the space of (isotopy classes of) decorated crowned surfaces. Specifically, it is defined as the space of pairs $(X,\textbf{h})$, where $X$ is a marked crowned hyperbolic metric on $\Sigma$ (i.e.: an element of $\mathcal{T}_\Sigma(\vec{0})$), and $\textbf{h}$ is a collection of horocycles on $X$ --- one for each cusp/tine.
\end{definition}

\begin{definition}[lambda length coordinates]
Consider an ideal triangulation $\triangle$ on $\Sigma$, for each arc $\alpha$ in $\triangle$, we define a lambda length
\[
\lambda_\alpha:
\widetilde{\mathcal{T}}_{\Sigma}
\rightarrow
\mathbb{R}_{>0}
\]
that assigns to a decorated marked hyperbolic surface $(X,\textbf{h})$ the lambda length of the lifts $\widetilde{\alpha}$ of the geodesic representative of $\alpha$ to the universal cover of $X$ in $\HH$, decorated by corresponding lifts of $\textbf{h}$.
\end{definition}

Penner introduced these spaces in \cite{pennercoords}, and showed that lambda length coordinates define global real analytically compatible coordinates for the decorated Teichm\"uller space of punctured (i.e.: cusped) surfaces. This also extends to crowned hyperbolic surfaces (see \cite[Theorem~2.2.25]{pennerbook}).

\section{Weil--Petersson volume form}
\label{sec:WPvol}

Teichm\"uller spaces $\mathcal{T}_{g,m}(\vec{b})$ of genus $g$ surfaces with $m$ cuffs of specified lengths are naturally symplectic via the \emph{Weil--Petersson symplectic form} \cite{weil}, with Darboux coordinates are given \emph{globally} (and somewhat miraculously) by Fenchel--Nielsen coordinates \cite{wolpertwpform}. The Weil--Petersson form is mapping class group invariant, and so its top exterior product defines the \emph{Weil--Petersson volume} on the moduli space $\mathcal{M}_{g,m}(\vec{b})$ --- this forms the backdrop for Mirzakhani's work on determining the Weil--Petersson volumes of $\mathcal{M}_{g,m}(\vec{b})$.\medskip

In contrast, Teichm\"uller spaces of crowned surfaces (either with or without fixed neck holonomy) are not naturally symplectic (they can even be odd-dimensional). Nevertheless, one can define a mapping class group invariant volume form on $\mathcal{T}_\Sigma(\vec{b})$ and $\mathcal{T}_\Sigma(\vec{b}|\vec{d})$ which generalise the Weil--Petersson volume form. Indeed, here are two approaches:

\begin{description}
    \item[Method~1] we can ``formally'' generalise classical expressions for the Weil--Petersson volume based on Fenchel--Nielsen parameters and shearing coordinates. This is certainly unsurprising to experts (see, e.g., Goncharov--Sun's \cite[Lemma~2.8]{goncharovsun}).
    
    \item[Method~2] Chekhov takes a more sophisticated approach in \cite{chekhov2024} using the Poisson algebraic structure of extended shearing coordinates on decorated Teichm\"uller space coming from the Goldman bracket, and exploiting the property, found in the classical setting, that the Weil--Petersson symplectic form is inverse to the Poisson bracket, to define the Weil--Petersson volume in terms of the Pfaffian of the matrix representing the Goldman bracket. 
\end{description}

In this section, we verify that the two methods yield the same volume form (up to sign), before clarifying the relationship between the Mirzakhani volumes of $\mathcal{M}_{\Sigma}(\vec{b}|\vec{d})$ and $\mathcal{M}_{\Sigma}(\vec{b})$. Finally, we verify Chekhov's assertion that that naively integrating the Weil--Petersson form over moduli spaces of crowned hyperbolic surfaces yields infinity precisely when one of the boundary of $\Sigma$ has more than one puncture.

\subsection{WP volume form on the moduli space of crowned surfaces}

We now define a generalization of the Weil--Petersson volume form (up to sign) on $\mathcal{M}_{\Sigma}(\vec{b}|\vec{d})$ and $\mathcal{M}_{\Sigma}(\vec{b})$ by first defining it for moduli spaces of $n$-gons and $n$-crowns, and then for general surfaces.

\subsubsection{The $n$-gon case}

\begin{definition}[Weil--Petersson volume form for $\mathcal{M}_{\mathbb{D}_n}$]
Let $\triangle$ be an arbitrary ideal triangulation of $\mathbb{D}_n$, let $\{s_\alpha\}_{\alpha\in\triangle}$ denote the associated shearing parameters for the shearing coordinates $s_\triangle:\mathcal{T}_{\mathbb{D}_n}\to\mathbb{R}^{\triangle}$, define the \emph{Weil--Petersson volume form} $\Omega^{\mathrm{WP}}_{\mathbb{D}_n}$ on $\mathcal{T}_{\mathbb{D}_n}$ as, up to sign, by
\begin{align}
\label{eq:vol-form-n-gon}  
\Omega^{\mathrm{WP}}_{\mathbb{D}_n} 
:= \pm \bigwedge_{\alpha\in\triangle}\mathrm{d}s_\alpha.
\end{align}
Since the (pure) mapping class group for $\mathbb{D}_n$ is trivial, this defines a volume form on $\mathcal{M}_{\mathbb{D}_n}=\mathcal{T}_{\mathbb{D}_n}$.
\end{definition}

\begin{remark}
The standard strategy to show that $\Omega^{\mathrm{WP}}_{\Sigma}$ is independent of the choice of the ideal triangulation is to verify that it is invariant under flips (i.e.: replace one edge by the opposite diagonal of the quadrilateral comprised of the two ideal triangles containing the given edge), and to show that you can deform from any ideal triangulation to another via a finite sequence of flips. We refer to \cite[Theorem~4.7]{pennerbook} for formulae related to the flipping process which can be used to show that $\Omega^{\mathrm{WP}}_{\Sigma}$ is flip-invariant up to sign.
\end{remark}

We now give a brute-force verification that the generalization of the Weil--Petersson volume form defined via \Cref{eq:vol-form-n-gon} is equal to that described in \cite[Equation~{(6.9)}]{chekhov2024}:

\begin{proposition}
\label{prop:vol-form-n-gon}
Consider points $\{z_j\}_{j=2,\ldots,n-2}$ in $(0,1)\subset\RR\cup\{\infty\}=\partial\HH$ parametrising an $n$-gon by specifying the positions of $n-3$ consecutive ideal vertices, where the $3$ other vertices are normalised so as to be placed at $0,1,\infty$. Then, $\{z_j\}_{j=2,\ldots,n-2}$ forms a coordinate chart on $\mathcal{M}_{\mathbb{D}_n}$, and
\begin{align}
\Omega^{\mathrm{WP}}_{\mathbb{D}_n}
=
\frac{\mathrm{d}z_2\wedge \mathrm{d}z_2\wedge\cdots\wedge \mathrm{d}z_{n-2}}
{z_2(z_3-z_2)(z_4-z_3)\cdots (z_{n-2}-z_{n-3})(1-z_{n-2})}.
\label{eq:chekhovngoncase}
\end{align}

\end{proposition}

\begin{proof}
Let $s_1, \cdots, s_{n-3}$ be the shearing coordinates obtained for the triangulation for $\mathbb{D}_n$ with all diagonals meeting at vertex $z_0 = z_n = \infty$, so that $s_i$ is assigned to the arc joining $\infty$ and $z_{i+1}$. Since \Cref{eq:vol-form-n-gon} holds for any triangulation of $\mathbb{D}_n$, we have
\[
\Omega^{\mathrm{WP}}_{\mathbb{D}_n}
= 
\mathrm{d}s_1 \wedge \cdots \wedge \mathrm{d}s_{n-3},
\]
up to sign. Introduce the following variables:
\[
y_2 = z_2,
\quad
y_i = z_{i}-z_{i-1},
\]
for $i = 3, \cdots, n-2$. Then 
\[
1-z_{n-2} = 1- \sum_{j=2}^{n-2}y_j
\quad\text{and}\quad 
\mathrm{d}y_2\wedge \mathrm{d}y_2\wedge\ldots\wedge \mathrm{d}y_{n-2}
= 
\mathrm{d}z_2\wedge \mathrm{d}z_2\wedge\ldots\wedge \mathrm{d}z_{n-2}.
\]
Hence \Cref{eq:chekhovngoncase} is equal to
\begin{align}
\label{eq: dOmega}
\frac{\mathop{\wedge}\limits_{i=2}^{n-2} \mathrm{d}y_i }{y_2y_3\cdots y_{n-2}(1- \sum_{j=2}^{n-2}y_j)}.
\end{align}
The shearing coordinates satisfy
\begin{align*}
s_1 &= \log \tfrac{z_3-z_2}{z_2} = \log \tfrac{y_3}{y_2},\\
s_j &= \log \tfrac{z_{j+2}-z_{j+1}}{z_{j+1}-z_j} = \log \tfrac{y_{j+2}}{y_{j+1}} \,\quad \text{for }j = 2, \cdots, n-4,\\
s_{n-3} &= \log \tfrac{1-z_{n-2}}{z_{n-2}-z_{n-3}} = \log \tfrac{1- \sum_{j=2}^{n-2}y_j}{y_{n-2}}.
\end{align*}
Thus, $\Omega^{\mathrm{WP}}_{\mathbb{D}_n}=\mathrm{d}s_1 \wedge \cdots \wedge \mathrm{d}s_{n-3}$ is given by
\begin{equation}
\label{eq: dOmega'}
\begin{split}
& (\mathrm{d}\log y_3 - \mathrm{d}\log y_2)
\wedge 
(\mathrm{d}\log y_4 - \mathrm{d}\log y_3) 
\wedge \cdots \wedge 
(\mathrm{d}\log y_{n-2} - \mathrm{d}\log y_{n-3}) \\ 
& \wedge \Bigl(\mathrm{d}\log \Bigl(1- \sum_{j=2}^{n-2}y_j \Bigr) - \mathrm{d}\log y_{n-2} \Bigr).
\end{split}
\end{equation}
Expand the brackets in \Cref{eq: dOmega'} and observe that the only non-zero wedge products are obtained as follows: for $k=0,\cdots, n-3$, 
\begin{itemize}
\item in the first $k$ brackets take the second summand, and
\item in the remaining $n-k-3$ brackets take the first summand.
\end{itemize}
The corresponding term equals $(-1)^{n-3}\mathop{\wedge}\limits_{j=2}^{n-2}\mathrm{d}\log y_j$ if $k=n-3$, and if $k<n-3$ the corresponding term is
\begin{align}
& \Bigl((-1)^k \mathop{\wedge}\limits_{i=2}^{k+1} \mathrm{d}\log y_i \Bigr) 
\wedge 
\Bigl(\mathop{\wedge}\limits_{i=k+3}^{n-2} \mathrm{d}\log y_i \Bigr)  
\wedge 
\Bigl(\mathrm{d}\log \Bigl(1- \sum_{j=2}^{n-2}y_j \Bigr)\Bigr) \\ 
&= 
(-1)^{k+(n-2-(k+3)+1)+1} \frac{y_{k+2} \cdot \mathop{\wedge}\limits_{i=2}^{n-2} \mathrm{d}\log y_i}{1- \sum_{j=2}^{n-2}y_j}
= 
(-1)^{n-3} \frac{y_{k+2} \cdot \mathop{\wedge}\limits_{i=2}^{n-2} \mathrm{d}\log y_i}{1- \sum_{j=2}^{n-2}y_j} .
\end{align}
Summing these $n-2$ terms, we obtain
\[
(-1)^{n-3} \mathop{\wedge}\limits_{i=2}^{n-2} 
\mathrm{d}\log y_i \Bigl(1+ \frac{y_2+y_3+\cdots+y_{n-2}}{1- \sum_{j=2}^{n-2}y_j}\Bigr) 
= (-1)^{n-3} \frac{\mathop{\wedge}\limits_{i=2}^{n-2} \mathrm{d}\log y_i}{1- \sum_{j=2}^{n-2}y_j},
\]
which equals \Cref{eq: dOmega}, up to sign.
\end{proof}

\subsubsection{The $n$-crown case}

\begin{definition}[Weil--Petersson volume form for $\mathcal{M}_{\mathbb{A}_n}$ and $\mathcal{M}_{\mathbb{A}_n}(d)$]
Let $\triangle$ be an arbitrary ideal triangulation of $\mathbb{A}_n$, let $\{s_\alpha\}_{\alpha\in\triangle}$ denote the associated shearing parameters for the shearing coordinates $s_\triangle:\mathcal{M}_{\mathbb{A}_n}(d)\to\mathbb{R}^{\triangle}$ (see \Cref{subsec:shears}), define the \emph{Weil--Petersson volume form} $\Omega^{\mathrm{WP}}_{\mathbb{A}_n}$ on $\mathcal{M}_{\mathbb{A}_n}=\mathcal{T}_{\mathbb{A}_n}$, up to sign, by
\begin{align}
\label{eq:vol-form-crowned-surf}  
\Omega^{\mathrm{WP}}_{\mathbb{A}_n} 
:= \pm \bigwedge_{\alpha\in\triangle}\mathrm{d}s_\alpha.
\end{align}
For $\mathcal{M}_{\mathbb{A}_n}(d)=\mathcal{T}_{\mathbb{A}_n}(d)$, let $\alpha_0\in\triangle$ be an arc incident to the interior puncture $p$ of $\mathbb{A}_n$. Then, we define
\begin{align}
\label{eq:vol-form-crowned-surf-neck}  
\Omega^{\mathrm{WP}}_{\mathbb{A}_n}(d) 
:= \pm \bigwedge_{\alpha\in\triangle\setminus\{\alpha_0\}}
\mathrm{d}s_\alpha.
\end{align}
Up to sign, this is independent of the choice of $\alpha_0$ because $d$ is the sum all of the $s_\alpha$ for which $\alpha\in\triangle$ is incident to $p$.
\end{definition}

To begin with, we do a first-principles check that the expression for $\Omega^{\mathrm{WP}}_{\mathbb{A}_n}(d)$ given here is equal to Chekhov's expression for the Weil--Petersson volume form on the moduli space of $n$-crowns with fixed neck length given in \cite[Lemma~{2.2}]{chekhov2024}:

\begin{proposition}[WP volume for $n$-crowns in shears]
\label{prop:vol-form-n-crown}
Consider $\{x_i\}_{i=1,\dotsc,n-1}$, where $x_i$ as the ratio of the lambda length (with respect to some decoration) of the arc from $q_{1}$ to $q_{i+1}$ to the lambda length of the arc from $q_{i+1}$ back to $q_{1}$, both taken in the counterclockwise direction \cite[Equation~{(2.6)}]{chekhov2024}. Then,
\begin{equation}
\label{eq:vol-form-crown-chekhov}
\Omega^{\mathrm{WP}}_{\mathbb{A}_n}(d)
=
\dfrac{ \mathrm{d}x_1\wedge\mathrm{d}x_2\wedge\ldots\wedge\mathrm{d}x_{n-1}}
{x_1(x_2-x_1)\cdots (x_i-x_{i-1})\cdots (x_{n-1}-x_{n-2}) }.
\end{equation}
Note that $\{x_i\}_{i=1,\dotsc,n-1}$ does not depend on the decorating horocycles.
\end{proposition}

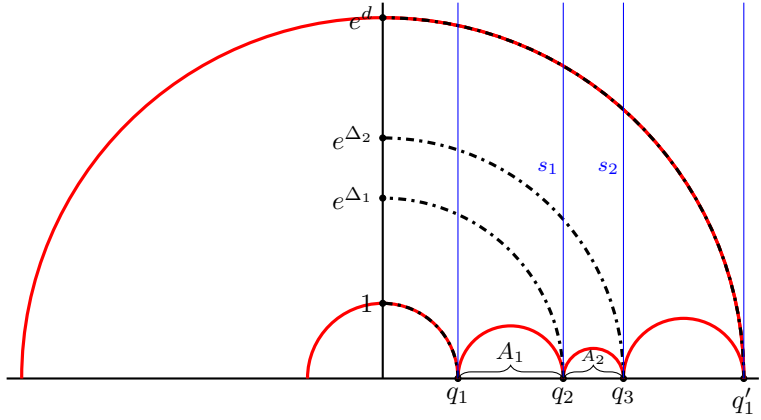
\begin{figure}[H]
    \centering
    \begin{tikzpicture}[scale=2]


  \draw[thick] (-2.5,0) -- (2.5,0);

  \draw[thick, black] (0,0) -- (0,2.5);

  \filldraw[black] (0.5,0) circle (0.02) node[below] {$q_1$};
  \filldraw[black] (1.2,0) circle (0.02) node[below] {$q_2$};
  \filldraw[black] (1.6,0) circle (0.02) node[below] {$q_3$};
  \filldraw[black] (2.4,0) circle (0.02) node[below] {$q_1'$};

  \draw[very thick, red] (1.2,0) arc (0:180:0.35);
  \draw[very thick, red] (1.6,0) arc (0:180:0.2);
  \draw[very thick, red] (2.4,0) arc (0:180:0.4);

  \draw[very thick, red] (0.5,0) arc (0:180:0.5);
  \draw[very thick, red] (2.4,0) arc (0:180:2.4);
  \draw[very thick, red] (2.4,0) arc (0:90:2.4);
  \draw[very thick, dash dot] (0.5,0) arc (0:90:0.5);
  \draw[very thick, dash dot] (1.2,0) arc (0:90:1.2);
  \draw[very thick, dash dot] (1.6,0) arc (0:90:1.6);
  \draw[very thick, dash dot] (2.4,0) arc (0:90:2.4);

  \filldraw[black] (0,0.5) circle (0.02) node[left] {$1$};
  \filldraw[black] (0,1.2) circle (0.02) node[left] {$e^{\Delta_1}$};
  \filldraw[black] (0,1.6) circle (0.02) node[left] {$e^{\Delta_2}$};
  \filldraw[black] (0,2.4) circle (0.02) node[left] {$e^{d}$};

  \draw[blue] (0.5,0) -- (0.5,2.5); 
  \draw[blue] (1.2,0) -- (1.2,2.5);
  \draw[blue] (1.6,0) -- (1.6,2.5);
  \draw[blue] (2.4,0) -- (2.4,2.5);

\draw [decorate,decoration={brace,amplitude=5pt}]
    (0.5,0) -- (1.2,0) node [midway,yshift=10pt] {\small{$A_1$}};  

\draw [decorate,decoration={brace,amplitude=5pt}]
    (1.2,0) -- (1.6,0) node [midway,yshift=8pt] {\tiny{$A_2$}}; 


\node [font=\fontsize{8}{8}, blue] at (1.1,1.4) {$s_1$};

\node [font=\fontsize{8}{8}, blue] at (1.5,1.4) {$s_2$};

\end{tikzpicture}
    \caption{A fundamental domain (enclosed by thick \textcolor{red}{red} lines) in the upper half-plane for a $3$-tined hyperbolic crown. In \textcolor{blue}{blue}: lifts of geodesic arcs that spiral around the neck curve.}
    \label{fig:crown-upper-half-plane}
\end{figure}

\begin{proof}
As in the proof of \Cref{prop:vol-form-n-gon}, we perform the calculation in the upper half-plane. Suppose that $d>0$. Following \cite[Page~7]{chekhov2024}, for each $i=1,\dotsc,n-1$ let $\Delta_i$ be the hyperbolic length of the counterclockwise segment of the geodesic neck connecting the nearest-point projections of $q_1$ and $q_{i+1}$. These satisfy $d>\Delta_{i-1}>\dotsc>\Delta_1>0$ (see \Cref{fig:crown-upper-half-plane} for an example with $n=3$).  We introduce the following variables:
\[
A_1 = e^{\Delta_1}-1, A_2 = e^{\Delta_2}-e^{\Delta_1}, \cdots, A_{n-1} = e^{\Delta_{n-1}}-e^{\Delta_{n-2}}.
\]
Further, let $Q = e^{d}-1$. Then 
\begin{equation}
\label{eq:Q-sumA}
e^{d}-e^{\Delta_{i}} = Q-\sum_{j=1}^i A_j, \, i=1,\cdots,n-1.
\end{equation}

Now let $s_1, \dotsc, s_n$ be the shearing coordinates associated with the ideal triangulation of $\AA_n$ by $n$ arcs emanating from the interior puncture, so that $s_i$ corresponds to the arc connecting the punctures $p_1$ and $q_{i+1}$. Then, we have (see \Cref{fig:crown-upper-half-plane}):
\begin{align*}
s_i &= \log \frac{A_{i+1}}{A_i},\,\quad \text{for }i = 1, \cdots, n-2,
\\
s_{n-1} &= \log\frac{Q-\sum_{j=1}^{n-1} A_j}{A_{n-1}}.
\end{align*}
Thus, $\mathrm{d}s_1 \wedge \cdots \wedge \mathrm{d}s_{n-1}$ is given by
\begin{equation}
\label{eq:vol-form-crown-A-var}   
\begin{split}
& (\mathrm{d}\log A_2-\mathrm{d}\log A_1)\wedge\cdots\wedge(\mathrm{d}\log A_{n-1}-\mathrm{d}\log A_{n-2})  \\ & \wedge\Bigl(\mathrm{d}\log \Bigl(Q-\sum_{j=1}^{n-1} A_j \Bigr)-\mathrm{d}\log A_{n-1}\Bigr).
\end{split}
\end{equation}
Notice that \Cref{eq:vol-form-crown-A-var} is almost the same as \Cref{eq: dOmega'} from the proof of \Cref{prop:vol-form-n-gon}. Then by a similar calculation, we obtain:
\begin{equation}
\label{eq:vol-form-crown-A-var-2} 
\mathrm{d}s_1 \wedge \cdots \wedge \mathrm{d}s_{n-1} = (-1)^{n-1}\frac{Q\cdot\mathop{\wedge}\limits_{i=1}^{n-1} \mathrm{d}A_i }{A_1\cdots A_{n-1}\Bigl(Q-\sum_{j=1}^{n-1} A_j \Bigr)},
\end{equation}
note the extra factor of $Q$ in the numerator in contrast to \Cref{eq: dOmega}.\medskip

Looking at \Cref{eq:vol-form-crown-chekhov}, observe that it can be rewritten as:
\begin{equation}
\label{eq:vol-form-crown-chekhov-2}
\Omega^{\mathrm{WP}}_{\mathbb{A}_n}(d)  = \mathrm{d}\log x_1 \wedge \mathrm{d}\log (x_2-x_2)\wedge\cdots\wedge \mathrm{d}\log(x_{n-1}-x_{n-2}).
\end{equation}
From \cite[Page~10]{chekhov2024}, Line $2$, we have
\[
x_i = e^{d/2} \frac{e^{\Delta_i}-1}{e^{d}-e^{\Delta_i}}, \,\quad \text{for }i = 1, \cdots, n-1.
\]
Then 
\[
x_{i+1}-x_i =
e^{d/2}\frac{(e^d-1)(e^{\Delta_{i+1}}-e^{\Delta_{i}})}{(e^{d}-e^{\Delta_{i+1}})(e^{d}-e^{\Delta_{i}})}, \,\quad \text{for }i = 1, \cdots, n-2.
\]
Since $d$ is a fixed constant, we have
\begin{align}
\label{eq:dlogx}
\mathrm{d}\log x_1 =\mathrm{d}\log\frac{e^{\Delta_1}-1}{e^{d}-e^{\Delta_1}}; \quad  \mathrm{d}\log (x_{i+1}-x_i) = \mathrm{d}\log \frac{e^{\Delta_{i+1}}-e^{\Delta_{i}}}{(e^{d}-e^{\Delta_{i+1}})(e^{d}-e^{\Delta_{i}})}.
\end{align}
We express \Cref{eq:dlogx} in the variables $A_1,\dotsc,A_{n-1}$ using \Cref{eq:Q-sumA}:
\[
\mathrm{d}\log x_1 = \mathrm{d}\log\frac{A_1}{Q-A_1}; \,\, \mathrm{d}\log (x_{i+1}-x_i) = \mathrm{d}\log \frac{A_{i+1}}{(Q-\sum_{j=1}^{i+1} A_j)(Q-\sum_{j=1}^i A_j)}.
\]
Observe that we have
\begin{equation}
\label{eq:dlog(x-x)}
\mathrm{d}\log (x_{i+1}-x_i) = \mathrm{d}\log \frac{A_{i+1}}{Q-\sum_{j=1}^{i+1} A_j} +\mathrm{d}\log \frac{1}{Q-\sum_{j=1}^{i} A_j},
\end{equation}
and the last term in \Cref{eq:dlog(x-x)} contributes trivially to the wedge product in \Cref{eq:vol-form-crown-chekhov-2} as it depends only on $A_1, \cdots, A_i$. Thus, \Cref{eq:vol-form-crown-chekhov-2} implies
\begin{align}
\label{eq:vol-form-crown-A-var-3}
\Omega^{\mathrm{WP}}_{\mathbb{A}_n}(d) = \mathrm{d}\log\frac{A_1}{Q-A_1} \wedge \mathrm{d}\log\frac{A_2}{Q-A_1-A_2} \wedge \cdots \wedge \mathrm{d}\log\frac{A_{n-1}}{Q-\sum_{j=1}^{n-1} A_j}. 
\end{align}
Next, observe that
\begin{align*}
\mathrm{d}\log \frac{A_i}{Q-\sum_{j=1}^{i} A_j} & = \mathrm{d}\log A_i - \mathrm{d}\log \Bigl(Q-\sum_{j=1}^{i} A_j\Bigr)= \frac{\mathrm{d}A_i}{A_i} + \frac{\mathrm{d}A_1+\cdots+\mathrm{d}A_i}{Q-\sum_{j=1}^{i} A_j}
\\ & = 
\mathrm{d}f(A_1,\cdots,A_{i-1})+  \left( \frac{1}{A_i}+ \frac{1}{Q-\sum_{j=1}^{i} A_j} \right)  \mathrm{d}A_i
\\ & = \mathrm{d}f(A_1,\cdots,A_{i-1})+ \frac{Q-\sum_{j=1}^{i-1} A_j}{Q-\sum_{j=1}^{i} A_j}\mathrm{d}\log A_i,
\end{align*}
where $\mathrm{d}f(A_1,\cdots,A_{i-1})$ is a $1$-form that depends only on $A_1,\cdots,A_{i-1}$ and thus contributes trivially to the wedge product in \Cref{eq:vol-form-crown-A-var-3}. Hence, we have
\[
\Omega^{\mathrm{WP}}_{\mathbb{A}_n}(d) = \prod_{i=1}^{n-1} \frac{Q-\sum_{j=1}^{i-1} A_j}{Q-\sum_{j=1}^{i} A_j} \bigwedge_{i=1}^{n-1} \mathrm{d}\log A_i  = \frac{Q}{Q-\sum_{j=1}^{n-1} A_j}\bigwedge_{i=1}^{n-1} \mathrm{d}\log A_i,
\]
which matches the expression in \Cref{eq:vol-form-crown-A-var-2}, up to sign.\medskip

If $d=0$, the quantities $\Delta_i$ vanish. To address this, we position a fundamental domain of the crown in the upper half-plane with $p_1 = \infty$ (as in \Cref{fig:crown-upper-half-plane}, in \textcolor{blue}{blue}). We then redefine $A_i, i=1,\dotsc,n-1$ to be Euclidean distance between $q_i$ and $q_{i+1}$ along $\partial\HH$, and let $Q$ be the Euclidean distance between $q_1$ and $q_1'.$ It is straightforward to verify that 
\[
x_i = \frac{\sum_{j=1}^i A_j}{Q-\sum_{j=1}^i A_j}, \,\quad \text{for }i = 1, \cdots, n-1.
\]
Hence,
\[
x_{i+1}-x_i = \frac{QA_{i+1}}{(Q-\sum_{j=1}^{i+1} A_j)(Q-\sum_{j=1}^{i} A_j)},
\]
and the preceding calculation for $d>0$ carries over with essentially no change.
\end{proof}

\subsubsection{All other cases (apart from $n$-gons)}

\begin{definition}[Weil--Petersson volume form for moduli spaces of general crowned hyperbolic surfaces]
    Given $\Sigma=\Sigma_{m,\vec{a}}$ of genus $g$, denote the components of $\Sigma$ its necks $\nu_1,\ldots,\nu_l$ to obtain the following decomposition 
\[
\Sigma\setminus(\nu_1\cup\ldots\cup\nu_l)=\Sigma_{g,m+l}\sqcup\mathbb{A}_{a_1}\sqcup\ldots\sqcup\mathbb{A}_{a_l}
\]
into a genus $g$ surface $\Sigma_{g,m+l}$ with $m+l$ interior punctures and a collection of crowns. Let $\{(\ell_k,\tau_k)_{k=1,\ldots,l}$ be respective Fenchel--Nielsen coordinates for the necks $\nu_1,\ldots,\nu_l$. Then, we define the Weil--Petersson volume form on $\mathcal{T}_{\Sigma}(\vec{b})$ and $\mathcal{T}_{\Sigma}(\vec{b}|\vec{d})$, up to sign, by
\begin{equation}
\begin{split}
\Omega^{\mathrm{WP}}_{\Sigma}(\vec{b})
& := \pm
\Omega^{\mathrm{WP}}_{\Sigma_{g,m+l}}(\vec{b},\vec{\ell})
\wedge 
\left(\bigwedge_{k=1}^l(\Omega^{\mathrm{WP}}_{\mathbb{A}_{a_1}}(\ell_k)
\wedge\mathrm{d}\ell_k\wedge\mathrm{d}\tau_k)\right),\\
\text{and}\quad
\Omega^{\mathrm{WP}}_{\Sigma}(\vec{b}|\vec{d})
& := \pm
\Omega^{\mathrm{WP}}_{\Sigma_{g,m+l}}(\vec{b},\vec{d})
\wedge  
\left(\bigwedge_{k=1}^l 
(\Omega^{\mathrm{WP}}_{\mathbb{A}_{a_1}}(d_1)\wedge\mathrm{d}\tau_k)
\right).
\label{prop:vol-form-scc-new}
\end{split}
\end{equation}
\end{definition}

\begin{proposition}
The Weil--Petersson volume forms $\Omega^{\mathrm{WP}}_{\Sigma}(\vec{b})$ and $\Omega^{\mathrm{WP}}_{\Sigma}(\vec{b}|\vec{d})$ are, up to sign, invariant with respect to the (pure) mapping class group.
\end{proposition}

\begin{proof}
    We prove this for $\Omega^{\mathrm{WP}}_{\Sigma}(\vec{b})$, the argument for $\Omega^{\mathrm{WP}}_{\Sigma}(\vec{b}|\vec{d})$ is essentially identical. The mapping class group $\mathrm{MCG}_\Sigma$ necessarily preserves the neck curves $\nu_1,\ldots,\nu_l$, and hence factors as 
    \begin{align*}
    \mathrm{MCG}_{\Sigma}
    &=
    \mathrm{MCG}_{\Sigma_{g,m+l}}
    \oplus
    \mathrm{MCG}_{\mathbb{A}_{a_1}}
    \oplus\cdots\oplus
    \mathrm{MCG}_{\mathbb{A}_{a_l}}
    \oplus
    \mathbb{Z}_{\nu_1}
    \oplus\ldots\oplus
    \mathbb{Z}_{\nu_l},\\
    &=
    \mathrm{MCG}_{\Sigma_{g,m+l}}
    \oplus
    \mathbb{Z}_{\nu_1}
    \oplus\ldots\oplus
    \mathbb{Z}_{\nu_l},
    \end{align*}
    where $\mathbb{Z}_{\nu_k}$ is the group of Dehn twists of $\nu_k$. Since the $\mathrm{MCG}_{\Sigma_{g,m+l}}$ preserves $\Omega^{\mathrm{WP}}_{\Sigma_{g,m+l}}(\vec{b},\vec{\ell})$ and leaves the other factors of $\Omega^{\mathrm{WP}}_{\Sigma}(\vec{b})$ untouched, it preserves $\Omega^{\mathrm{WP}}_{\Sigma}(\vec{b})$. Similarly, the $\mathbb{Z}_{\nu_k}$ preserves $\mathrm{d}\tau_k$ and leaves all other terms untouched and hence also preserves $\Omega^{\mathrm{WP}}_{\Sigma}(\vec{b})$. Since they generate $\mathrm{MCG}_{\Sigma}$, this means that $\Omega^{\mathrm{WP}}_{\Sigma}(\vec{b})$ is invariant under the action of $\mathrm{MCG}_{\Sigma}$.
\end{proof}

\begin{remark}
    Although Chekhov does not state explicitly in \cite{chekhov2024} a statement equivalent to \Cref{prop:vol-form-scc-new}, it is implicitly used as a part of deriving \cite[Lemma~2.6]{chekhov2024}. 
\end{remark}

\section{Mirzakhani volumes of moduli spaces of crowned hyperbolic surfaces with neck length constraints}
\label{sec:vol-mod-space-crowns}

In this section, we
\begin{itemize}
    \item carry out the calculation of the Mirzakhani volumes of $\mathcal{M}_{\mathbb{A}_n}(d)$ for all $n\geqslant1$ (\Cref{thm:fixedneckvol});
    \item determine the generating function for the Mirzakhani volumes of the moduli space of $n$-crowns (\Cref{prop:genfunction-crowns});
    \item determine the Mirzakhani volumes of $\mathcal{M}_{\Sigma}(\vec{b}|\vec{d})$ for all $\Sigma$ which are not $n$-gons (\Cref{cor:fixedneckvolumegeneral}).
\end{itemize}  
There are discrepancies between the answers we obtain and those obtained by Chekhov in \cite{chekhov2024}, and we account for the differences in \Cref{rmk:chekhoverror}.

\subsection{Chekhov action and its necessity}

Chekhov asserts on page~2 of \cite{chekhov2024} the Weil--Petersson volume form introduced in \Cref{sec:WPvol} can have infinite volume. Specifically,

\begin{proposition}

The integral 
\[
\int_{\mathcal{M}_{\Sigma}(\vec{b}|\vec{d})}
\Omega^{\mathrm{WP}}_{\Sigma}(\vec{b}|\vec{d})
\]
is finite if and only if $a_1=a_2=\ldots=a_l=1$, meaning that all other volumes of moduli spaces for crowned surfaces (of fixed neck length) are infinite.   
\end{proposition}

\begin{proof}
    \Cref{prop:vol-form-scc-new} implies:
    \[
    \int_{\mathcal{M}_{\Sigma}(\vec{b}|\vec{d})}
    \Omega^{\mathrm{WP}}_{\Sigma}(\vec{b}|\vec{d})
    =
    V_{g+l}(\vec{b},\vec{d})
    \prod_{k=1}^l
    \left(
    d_k\cdot\int_{\mathcal{M}_{\mathbb{A}_{a_k}}(d_k)}\Omega^{\mathrm{WP}}_{\mathbb{A}_{a_k}}(d_k).
    \right)
    \]
Invoking \Cref{eq:vol-form-crowned-surf-neck}, we see that $\int_{\mathcal{M}_{\mathbb{A}_{a_k}}(d_k)}\Omega^{\mathrm{WP}}_{\mathbb{A}_{a_k}}$ is finite if and only if $a_k=1$, and the result follows.
\end{proof}

Chekhov uses the infinitude of the naive Weil--Petersson integral as one reason for defining a renormalising action: 
\begin{definition}[Chekhov's action]
Given $X\in\mathcal{M}_{\Sigma}(\vec{b})$, choose an arbitrary horocyclic decoration on $X$ and let $h_j$ be the length of the horocyclic arc at $j$-th tine $q_j$ and $\lambda_j$ is the lambda length of the $j$-th arch $\alpha_j$. Define \emph{Chekhov's action} $S:\mathcal{M}_{\Sigma}(\vec{b})\rightarrow\mathbb{R}$ by
\[
S(X) := \log \left(\prod_{i=1}^n \lambda_{i}h_i \right).
\]
\end{definition}
Although $S$ is defined using horocyclic decorations, Chekhov shows \textit{a posteriori} that $S$ is decoration independent, and properly defines an action on $\mathcal{M}_{\Sigma}(\vec{b})$.

\begin{definition}[Mirzakhani volume]\label{defn:mirzvolume}
We refer to the normalised integral
\begin{align*}
    V_{\Sigma}(\vec{b}|\vec{d}):=
    \int_{\mathcal{M}_{\Sigma}(\vec{b}|\vec{d})} e^{-S(X)}\;\Omega^{\mathrm{WP}}_{\Sigma}(\vec{b}|\vec{d})
\end{align*}
as the \emph{Mirzakhani volume} of the moduli space $\mathcal{M}_{\Sigma}(\vec{b}|\vec{d})$ of crowned hyperbolic surfaces homeomorphic to $\Sigma$ with fixed cuff and neck lengths. Likewise, we define
\begin{align*}
    V_{\Sigma}(\vec{b})
    :=\int_{\mathcal{M}_{\Sigma}(\vec{b})} e^{-S(X)}\;\Omega^{\mathrm{WP}}_{\Sigma}(\vec{b})
\end{align*}
as the \emph{Mirzakhani volume} of the moduli space $\mathcal{M}_{\Sigma}(\vec{b})$.
\end{definition}

\begin{remark}
    One might (very reasonably) wonder where the name \emph{Mirzakhani volume} comes from, and whilst Chekhov does use this phrase in \cite[Section~2.3]{chekhov2024} entitled ``Mirzakhani volumes for $\mathcal{M}_{g,s,\mathbf{n}}$'' and expressions like ``the original Mirzakhani volumes'' (implying that there are new ones), there is some ambiguity as to whether he referred to the induced volume for his action as ``Mirzakhani volumes''. We choose to adopt this term both to distinguish it from the Weil--Petersson volume and in memory of Mirzakhani's work.
\end{remark}

\subsection{Mirzakhani volumes of moduli spaces of crowned hyperbolic surfaces with fixed neck lengths}

In \cite[Theorem~2.4]{chekhov2024}, Chekhov obtains an expression for the Mirzakhani volume of $\mathcal{M}_{\mathbb{A}_n}(d)$ as an integral over an $(n-1)$-dimensional simplex and uses this to assert the finiteness of these integrals. The coordinate system he utilizes comes from the distances $\Delta_i$ between the orthogonal projections of the punctures to the crown's neck. In particular, he explicitly carries out the calculation for the moduli space of $3$-crowns with fixed neck length \cite[Section~2.5]{chekhov2024}.\medskip 

We determine these volumes for all $n\geq1$, and hence determine 

\begin{theorem}
\label{thm:fixedneckvol}
The Mirzakhani volume $V_{\mathbb{A}_n}(d)$ is:
\begin{align*}
V_{\mathbb{A}_n}(d) 
=
\left\{
\begin{array}{rl}
\frac{d}{\sinh{d/2}} \cdot 
\frac{\prod_{j=1}^{k-1} (d^2+(2j)^2\pi^2)}{2(n-1)!},
&\quad\text{if }n=2k\text{ is even;}
\\[0.25em]
\frac{1}{\cosh{d/2}} \cdot  \frac{\prod_{j=1}^{k} (d^2+(2j-1)^2\pi^2 )}{2(n-1)!},
&\quad\text{if }n = 2k+1\text{ is odd.}
\end{array}
\right.
\end{align*}

\end{theorem}

\begin{proof}

To determine the volume of the moduli space for $n$-crowns $\mathcal{M}_{\mathbb{A}_{n}}(\vec{d})$, we need
\begin{itemize}
    \item a domain of integration,
    \item an expression for the Weil--Petersson volume form,
    \item and an expression for the Chekhov action.
\end{itemize}
Since the (pure) mapping class group for a boundary punctured disk with one interior puncture $\AA_n$ is trivial, the moduli space for $n$-crowns $\mathcal{M}_{\mathbb{A}_{n}}(\vec{d})$ is equal to the Teichm\"uller space $\mathcal{T}_{\mathbb{A}_{n}}(\vec{d})$. Fix the shearing coordinates $\{s_1,\dotsc,s_n\}$ on $\mathcal{T}_{\mathbb{A}_{n}}(\vec{d})$ as in the proof of \Cref{prop:vol-form-n-crown}. Then by \Cref{prop:shearingcoords}, the domain of integration is given by the hyperplane $H_d =\left\{\sum_{i=1}^n s_i =d\right\} \subset \RR^n$. The shearing coordinates based presentation of the Weil--Petersson volume form is given by
\[
\Omega^{\mathrm{WP}}_{\mathbb{A}_n}(d)= \mathrm{d}s_1 \wedge \cdots \wedge \mathrm{d}s_{n-1}.
\]
And the expression for the Chekhov action in shearing coordinates is given by Chekhov in \cite[Equation~5.1]{chekhov2024}:
\[
e^{-S} = \frac{1}{\prod_{i=1}^n 2\cosh (s_i/2)}.
\]
Then we can set up the volume integral:
\begin{align}
V_{\mathbb{A}_n}(d) 
= \int_{\mathcal{T}_{\mathbb{A}_{n}}(\vec{d})} e^{-S} \,\Omega^{\mathrm{WP}}_{\mathbb{A}_n}(d) 
= \int_{H_d} \,\, \frac{\mathrm{d}s_1 \wedge \cdots \wedge \mathrm{d}s_{n-1}}{ \prod_{i=1}^n 2\cosh (s_i/2)}.\label{eq:ncrownvolume}
\end{align}
This is a convolution of functions, specifically:
\begin{align}
V_{\mathbb{A}_n}(d) 
= \left(\underbrace{\frac{1}{2\cosh(s/2)} * \cdots * \frac{1}{2\cosh(s/2)}}_{n \text{ times}}\right)(d),
\end{align}
and we may apply the convolution theorem:
\begin{align}
\underbrace{f * f *  \cdots * f}_{{n \text{ times}}} = \mathcal{F}^{-1} (\mathcal{F}(f)^n),
\end{align}
where $f(x) = \frac{1}{2}\text{sech}(x/2)$ and $\mathcal{F}$ denotes the Fourier transform. A straightforward residue-based calculation (see also \cite[Equation~3.981(3)]{MR2360010}) yields:
\begin{align}
\mathcal{F}[f](\xi)  
= \int_\RR f(x)\ e^{-2\pi i \xi x}\; \mathrm{d}x 
= \frac{\pi}{\cosh{2\pi^2\xi}}. 
\end{align}
Thus, by the convolution theorem:
\[
V_{\mathbb{A}_n}(x) = 
\int_{\RR} \frac{\pi^n e^{2\pi i\xi x}}{\cosh^n{2\pi^2 \xi}} \mathrm{d}\xi .
\]
Let $\xi = \frac{\log t}{2\pi^2},$ then $\mathrm{d}\xi = \frac{\mathrm{d}t}{2\pi^2 t}$, and
\begin{equation*}
\begin{split}
V_{\mathbb{A}_n}(x) &  = \pi^n \int_{0}^\infty \frac{e^{2\pi i  \frac{\log t}{2\pi^2} x }}{\cosh^n({\log t})} \frac{\mathrm{d}t}{2\pi^2 t} = \frac{\pi^{n-2}}{2} \int_0^\infty \frac{t^{\frac{ix}{\pi}-1}}{\left(\frac{t+t^{-1}}{2}\right)^n}\mathrm{d}t  \\&= 2^{n-1}\pi^{n-2} \int_0^\infty \frac{t^{\frac{ix}{\pi}+n-1}}{(t^2+1)^n}\mathrm{d}t.
\end{split}
\end{equation*}
Let $u = \frac{1}{t^2+1}$. Then $t = (u^{-1}-1)^{1/2}$, and 
\[
\mathrm{d}t = \frac{1}{2}(u^{-1}-1)^{-1/2}\cdot (-1)\cdot u^{-2}\mathrm{d}u.
\] 
Therefore,
\begin{align*}
V_{\mathbb{A}_n}(x) & = (2\pi)^{n-2} \int_{0}^1 (u^{-1}-1)^{\frac{ix}{2\pi}+\frac{n}{2}-\frac{1}{2}}\cdot u^n \cdot  (u^{-1}-1)^{-1/2} u^{-2}\mathrm{d}u \\&
= (2\pi)^{n-2} \int_0^1 (u^{-1}-1)^{\frac{ix}{2\pi}+\frac{n}{2}-1} u^{n-2} \mathrm{d}u.
\\& = (2\pi)^{n-2} \int_0^1 (1-u)^{\frac{ix}{2\pi}+\frac{n}{2}-1}\cdot u^{-\frac{ix}{2\pi}+\frac{n}{2}-1}\mathrm{d}u \\&
= (2\pi)^{n-2} B \Bigl(-\frac{ix}{2\pi}+\frac{n}{2}, \frac{ix}{2\pi}+\frac{n}{2} \Bigr) \\& = \frac{(2\pi)^{n-2}}{\Gamma(n)}\Gamma\Bigl(-\frac{ix}{2\pi}+\frac{n}{2}\Bigr)\Gamma\Bigl(\frac{ix}{2\pi}+\frac{n}{2}\Bigr) \\& = \frac{(2\pi)^{n-2}}{\Gamma(n)} \Big|\Gamma\Bigl(\frac{n}{2}+\frac{ix}{2\pi}\Bigr)\Big|^2,
\end{align*}
where $B(z_1,z_1)$ is the beta function and $\Gamma(z)$ is the gamma function.\medskip

To finish up, we use the absolute value formulas for arguments of integer or half-integer real part \cite[Equation~1.2(8)-(9)]{batemanHTF1}.
Let $x = d$. Then for $n = 2k$:
\[
\Big|\Gamma\Bigl(\frac{n}{2}+\frac{id}{2\pi}\Bigr)\Big|^2 = \frac{d/2}{\sinh{d/2}} \prod_{j=1}^{k-1} \Bigl((d/2\pi)^2+j^2\Bigr).
\]
Hence
\[
\frac{(2\pi)^{n-2}}{\Gamma(n)} \Big|\Gamma\Bigl(\frac{n}{2}+\frac{id}{2\pi}\Bigr)\Big|^2 = \frac{d}{\sinh{d/2}} \cdot \frac{\prod_{j=1}^{k-1} \big(d^2+(2j)^2\pi^2\big)}{2(n-1)!}.
\]
For $n = 2k+1$, we have:
\[
\Big|\Gamma\Bigl(\frac{n}{2}+\frac{id}{2\pi}\Bigr)\Big|^2 = \frac{\pi}{\cosh{d/2}} \prod_{j=1}^{k} \bigl((d/2\pi)^2+(j-1/2)^2\bigr).
\]
Hence
\[
\frac{(2\pi)^{n-2}}{\Gamma(n)} \Big|\Gamma\Bigl(\frac{n}{2}+\frac{id}{2\pi}\Bigr)\Big|^2 = \frac{1}{\cosh{d/2}} \cdot \frac{\prod_{j=1}^{k} \bigl(d^2+(2j-1)^2\pi^2\bigr)}{2(n-1)!}.
\]
\end{proof}

Applying \cite[Lemma~2.6]{chekhov2024} (which follows easily from \Cref{prop:vol-form-scc-new}), we obtain:
\begin{corollary}
\label{cor:fixedneckvolumegeneral}
For $\Sigma=\mathbb{A}_{a_1,a_2}(d)$, the Mirzakhani volume of $\mathcal{M}_{\mathbb{A}_{a_1,a_2}}(d)$ is
    \[
    V_{\mathbb{A}_{a_1,a_2}}(d)
    =
    d\cdot V_{\mathbb{A}_{a_1}}(d)\cdot V_{\mathbb{A}_{a_2}}(d).
    \]
    For $\Sigma\neq\mathbb{D}_n,\mathbb{A}_n,\mathbb{A}_{a_1,a_2}$, let $g$ denote the genus of $\Sigma$. The volume of $\mathcal{M}_\Sigma(\vec{b}|\vec{d})$ is
    \[
    V_\Sigma(\vec{b}|\vec{d})
    =
    V_{g,m+l}(\vec{b},\vec{d})
    \cdot
    \prod_{k=1}^{l}d_k \cdot V_{\mathbb{A}_{a_k}}(d_k),
    \]
    where $V_{g,m+l}(\vec{b},\vec{d})\in\mathbb{Q}_{>0}[\pi^2,b_1^2,\ldots,b_m^2,d_1^2,\ldots,d_l^2]$ is the Weil--Petersson volume of the moduli space $\mathcal{M}_{g,m+l}(\vec{b},\vec{d})$.
\end{corollary}

\begin{remark}\label{rmk:chekhoverror}
The result of \Cref{thm:fixedneckvol} disagrees with Chekhov's calculation for $n = 3$ \cite[Example~2.5.1]{chekhov2024} (and with the case $n=2$ claimed in \cite[Lemma~2.6]{chekhov2024}; see \Cref{app:2-crowns} for a computation of $V_{\AA_2}(d)$ using $\lambda$-lengths). This discrepancy can be traced to a subtle computational oversight in \cite[Equation~2.4]{chekhov2024}: following the notation of \cite{chekhov2024}, for $i=0$ it should read
\begin{align*}
s_0 &= \sqrt{e^P r_0/r_{n-1}} + \sqrt{r_0/r_1} = \frac{2\sqrt{e^{P}r_0r_1}+2\sqrt{r_0r_{n-1}}}{2\sqrt{r_1r_{n-1}}}
\\&= \frac{e^{P/2}(e^{\Delta_1}-1)+e^{-P/2}(e^P-e^{\Delta_{n-1}})}{2\sqrt{r_1r_{n-1}}}
=\frac{e^{-P/2}(e^{P+\Delta_1}-e^{\Delta_{n-1}})}{2\sqrt{r_1r_{n-1}}}.
\end{align*}
Taking the product of horocycle lengths, we obtain
\begin{align*}
\prod_{j=0}^{n-1} s_j = \frac{\left(\prod_{i=1}^{n-1}(e^{\Delta_{i+1}}-e^{\Delta_{i-1}})\right)\cdot e^{-P/2}(e^{P+\Delta_1}-e^{\Delta_{n-1}}) }{\prod_{i=1}^n(e^{\Delta_{i}}-e^{\Delta_{i-1}})}  = e^{-P/2} \prod_{i=1}^n\frac{e^{\delta_i+\delta_{i+1}}-1}{e^{\delta_i}-1}.
\end{align*}
Taking the reciprocal, we obtain that the action functional $e^{-S}$ in \cite[Theorem~2.4]{chekhov2024} is missing a multiplicative factor of $e^{P/2}$. Moreover, the multiplicative factor $\frac{1}{2^{n-1}}$ in \cite[Theorem~2.4]{chekhov2024} is redundant: it cancels out upon rewriting \cite[Equation~2.9]{chekhov2024} from being a product of hyperbolic sine based terms to being expressed in terms of exponentials. Hence, the correct form of \cite[Theorem 2.4]{chekhov2024} should be
\begin{align*}
V^{\text{crown}}_{n,P}= \int_{\delta_1+\cdots +\delta_{n-1}\le P} \prod_{i=1}^{n-1}d\delta_i \frac{e^{P/2}(e^P-1)}{\prod_{i=1}^n \bigl( e^{\delta_i+\delta_{i+1}}-1 \bigr)}
\end{align*}
with $\delta_i>0$, $\delta_{i+n}=\delta_i$, and $\delta_1+\cdots +\delta_{n-1}+\delta_n=P$. 

\end{remark}

\subsection{Generating function}
We set the generating function to take the form
\[
P_{\mathbb{A}_{\bullet}}(x,d):=\sum_{n=1}^\infty V_{\mathbb{A}_n}(d)\,x^{n}.
\]

\begin{theorem}[generating function for volumes of moduli spaces of crowns]
\label{prop:genfunction-crowns}
The generating function $P_{\mathbb{A}_{\bullet}}(x,d)$ for the volumes of crowns is
\begin{align}
P_{\mathbb{A}_{\bullet}}(x,d)
= 
\frac{x}{\sqrt{1-\pi^2x^2}} \cdot\frac{\sinh(\frac{d}{2}+\frac{d}{\pi}\arcsin\pi x)}{\sinh d}.
\end{align}
\end{theorem}
\begin{proof}
We write $P_{\mathbb{A}_{\bullet}}(x,d) = P_{\mathbb{A}_{\bullet}}^{even}(x,d)+P_{\mathbb{A}_{\bullet}}^{odd}(x,d)$, where $P_{\mathbb{A}_{\bullet}}^{even}$ and $P_{\mathbb{A}_{\bullet}}^{odd}$ respectively correspond to even $n$ volume terms and odd $n$ volume terms, and we compute them separately. By \cite[Equation~1.414(1)]{MR2360010}: 
\begin{align}
\label{eq:cos(arcsinh)}
\cos (d \,\text{arcsinh}\, {x}) = 1 - \sum_{k=0}^\infty (-1)^k \frac{d^2(d^2+2^2)\dotsc(d^2+(2k)^2)}{(2k+2)!} x^{2k+2}.
\end{align}
Since $\text{arcsinh}(ix)=i \arcsin(x)$ and $\cos  ix = \cosh x$, we can express \Cref{eq:cos(arcsinh)} as
\begin{align}
\cosh (d \arcsin {x}) = 1 + \sum_{k=0}^\infty  \frac{d^2(d^2+2^2)\dotsc(d^2+(2k)^2)}{(2k+2)!} x^{2k+2}.
\end{align}
Together with \Cref{thm:fixedneckvol}, it follows that
\begin{equation}
\label{eq:G}
\begin{split}
\frac{\cosh \left(\frac{d}{\pi} \arcsin {\pi x}\right) -1}{2d\sinh(d/2)} &= \sum_{k=0}^\infty  \frac{d}{\sinh(d/2)}\frac{(d^2+2^2\pi^2)\dotsc(d^2+(2k)^2\pi^2)}{2(2k+2)!} x^{2k+2} 
\\&
= \sum_{k=0}^\infty  \frac{V_{\mathbb{A}_{2k+2}}(d)}{2k+2} x^{2k+2}.
\end{split}
\end{equation}
Let $G(x,d)= \frac{\cosh \left(\frac{d}{\pi} \arcsin {\pi x}\right) -1}{2d\sinh(d/2)}$. It follows from \Cref{eq:G} that
\begin{align}
P_{\mathbb{A}_{\bullet}}^{even}(x,d)=x \frac{\partial G}{\partial x},  
\end{align}
and we compute
\begin{align}
\label{eq: P even}
P_{\mathbb{A}_{\bullet}}^{even}(x,d) = \frac{x}{2\sqrt{1-\pi^2x^2}} \frac{\sinh\left(\frac{d}{\pi} \arcsin {\pi x}\right)}{\sinh(d/2)}.   
\end{align}
Next, by \cite[Equation~1.414(2)]{MR2360010}: 
\begin{align}
\label{eq:sin(arcsinh)}
\sin (d \,\text{arcsinh}\, {x}) = dx - d\sum_{k=1}^\infty (-1)^{k+1} \frac{(d^2+1^2)(d^2+3^2)\dotsc(d^2+(2k-1)^2)}{(2k+1)!} x^{2k+1}.
\end{align}
Since $\text{arcsinh}(ix)=i \arcsin(x)$ and $\sin ix= i\sinh x $, we can express  \Cref{eq:sin(arcsinh)} as

\begin{align}
\sin (d \,\text{arcsinh}\, {x}) = dx - d\sum_{k=1}^\infty (-1)^{k+1} \frac{(d^2+1^2)(d^2+3^2)\dotsc(d^2+(2k-1)^2)}{(2k+1)!} x^{2k+1}.
\end{align}
Together with \Cref{thm:fixedneckvol}, it follows that
\[
\begin{aligned}
&\frac{\sinh \left(\frac{d}{\pi} \arcsin {\pi x}\right) }{2d\cosh(d/2)} \\
\smash{\mathllap{=}}\, &\frac{x}{2\cosh(d/2)} + \sum_{k=1}^\infty  \frac{1}{\cosh(d/2)}\frac{(d^2+1^2\pi^2)(d^2+3^2\pi^2)\dotsc(d^2+(2k-1)^2\pi^2)}{2(2k+1)!} x^{2k+1} \\
\smash{\mathllap{=}}\, &   \sum_{k=0}^\infty  
\frac{V_{\mathbb{A}_{2k+1}}(d)}{2k+1} x^{2k+1}.
\end{aligned}
\]
Let $H(x,d) = \frac{\sinh \left(\frac{d}{\pi} \arcsin {\pi x}\right) }{2d\cosh(d/2)}$. Then  $P_{\mathbb{A}_{\bullet}}^{odd}(x,d) = x \frac{\partial H}{\partial x} $, and we get:

\begin{align}
\label{eq: P odd}
P_{\mathbb{A}_{\bullet}}^{odd}(x,d) = \frac{x}{2\sqrt{1-\pi^2x^2}}  \frac{\cosh\left(\frac{d}{\pi} \arcsin {\pi x}\right)}{\cosh(d/2)}.    
\end{align}
Combining \Cref{eq: P even} with \Cref{eq: P odd}, we conclude
\begin{equation}
\begin{split}
P_{\mathbb{A}_{\bullet}}(x,d) &= \frac{x}{2\sqrt{1-\pi^2x^2}} \left( \frac{\sinh\left(\frac{d}{\pi} \arcsin {\pi x}\right)}{\sinh(d/2)}+\frac{\cosh\left(\frac{d}{\pi} \arcsin {\pi x}\right)}{\cosh(d/2)}\right) 
\\&
= \frac{x}{\sqrt{1-\pi^2x^2}} \cdot\frac{\sinh(\frac{d}{2}+\frac{d}{\pi}\arcsin\pi x)}{\sinh d}.   
\end{split}
\end{equation}

\end{proof}

\section{Mirzakhani volumes of moduli spaces of crowned hyperbolic surfaces without neck length constraints}

The Mirzakhani volumes for moduli spaces of crowned hyperbolic surfaces of fixed neck lengths $d_1,\ldots,d_l$ take the form of something polynomial in $d_1,\ldots,d_l$ over something exponential, and thus the Mirzakhani volumes for moduli spaces without neck constraints are convergent. In this section, we
\begin{itemize}
    \item compute the Mirzakhani volumes of $\mathcal{M}_{\mathbb{A}_n}$ for all $n\geqslant1$ (\Cref{thm:crownfullvol});
    \item show that the Mirzakhani volumes of the moduli space of $(a_1,a_2)$-crowns is a rational polynomial in $\log2$ and  integral values of the Riemann zeta function (\Cref{thm:vol-biannuli});
    \item show that the Mirzakhani volumes of $\mathcal{M}_{\Sigma}(\vec{b})$ are rational polynomials in integral values of the Riemann zeta function and the Dirichlet beta function for all remaining $\Sigma$ which are not $n$-gons (\Cref{thm:generalmirzvolumes}).
\end{itemize}

\subsection{Mirzakhani volumes of moduli spaces of $n$-crowns} 

We first find the Mirzakhani volumes for full moduli spaces of $n$-crowns (i.e.: without neck length constraints):
\begin{proposition}
\label{thm:crownfullvol}
For $n \geqslant  1$,
\[
V_{\mathbb{A}_{n}}= \frac{\pi^n}{2}.
\]
\end{proposition}

\begin{proof}
Note that the function $f(x) = \frac{1}{2}\text{sech}(x/2)$ is even. Then \cite[Equation~3.511(1)]{MR2360010}) tells us:
\begin{align*}
V_{\mathbb{A}_{n}} &= \int_{\{\sum_{i=1}^n s_i\geqslant0\}} \bigwedge_{i=1}^n\frac{\mathrm{d}s_i}{ 2\cosh (s_i/2)} \\&= \frac{1}{2}\int_{\RR^n} \bigwedge_{i=1}^n\frac{\mathrm{d}s_i}{ 2\cosh (s_i/2)} = \frac{1}{2}\left(\int_\RR\frac{\mathrm{d}s}{2\cosh(s/2)}\right)^n = \frac{\pi^n}{2}. 
\end{align*}    
\end{proof}

\subsection{Mirzakhani volumes of moduli spaces of $(a_1,a_2)$-annuli}
\label{sec:vol-annuli}

\begin{theorem}
\label{thm:vol-biannuli}
If $n=a_1+a_2=\dim_{\mathbb{R}}\mathcal{M}_{\mathbb{A}_{a_1,a_2}}$ is odd, then
\[
V_{\mathbb{A}_{a_1,a_2}} =\sum_{i=0}^{\frac{n-3}{2}} \rho_{2i}(a_1,a_2) \pi^{2i}\zeta(n-2i)
\]
where $\rho_{2i}(a_1,a_2) \in \QQ_{>0}$ and $\zeta(s)$ is the Riemann zeta function. If $n=a_1+a_2$ is even, then
\[
V_{\mathbb{A}_{a_1,a_2}} = \sum_{i=0}^{\frac{n-4}{2}} \rho_{2i}(a_1,a_2) \pi^{2i}\zeta(n-1-2i)+ \rho_{n-2}(a_1,a_2)\pi^{n-2}\log2,
\]
where $\rho_{2i}(a_1,a_2) \in \QQ_{>0}$, and $\rho_{n-2}(a_1,a_2) \in \QQ_{\geqslant 0}$. The coefficient $\rho_{n-2}(a_1,a_2)$ is strictly positive if and only if both $a_1,a_2$ are odd.
\end{theorem}

\begin{proof}
By \cite[Lemma~2.6]{chekhov2024} (or just \Cref{prop:vol-form-scc-new}):
\begin{align}
\label{eq:int-vol-bi-annulus}
V_{\mathbb{A}_{a_1,a_2}} 
= \int_0^\infty\int_0^\ell 
V_{\mathbb{A}_{a_1}}(\ell)
\cdot V_{\mathbb{A}_{a_2}}(\ell)
\,\mathrm{d}\tau \mathrm{d}\ell 
=\int_0^\infty 
\ell
\cdot V_{\mathbb{A}_{a_1}}(\ell)
\cdot V_{\mathbb{A}_{a_2}}(\ell) 
\,\mathrm{d}\ell.
\end{align}

The specific integrals going forth depend on the parity of $a_1,a_2$, and we split into the following cases: 
\begin{description}
    \item[Case~1] $a_1,a_2$ have different parity, or equivalently, $n$ is odd. 
    \item[Case~2] $a_1,a_2$ are both even.
    \item[Case~3] $a_1,a_2$ are both odd.
\end{description}
By \Cref{thm:fixedneckvol}, these three cases respectively correspond to integrals of the form:
\begin{align*}
&\text{For Case~1,}\quad
&\int_0^\infty
\frac{P(\ell)\;\mathrm{d}\ell}{\sinh(\ell/2)\cosh(\ell/2)},&
\quad\text{$P$ is a degree $n-1$ even polynomial.}\\
&\text{For Case~2,}\quad
&\int_0^\infty
\frac{Q(\ell)\;\mathrm{d}\ell}{\sinh^2 (\ell/2)},&
\quad\text{$Q$ is a degree $n-1$ odd polynomial.}\\
&\text{For Case~3,}
&\int_0^\infty
\frac{R(\ell)\;\mathrm{d}\ell}
{\cosh^2 (\ell/2)},&
\quad\text{$R$ is a degree $n-1$ odd polynomial}.
\end{align*}
Moreover, by \Cref{thm:fixedneckvol}, if we regard $P,Q,R$ as polynomials $\mathbb{Q}[\pi^2,\ell]$ where $\pi^2$ has degree $2$, then these three polynomials are homogeneous with strictly positive rational coefficients. To finish the proof, it suffices to show that the constituent monomials in $P,Q,R$, when divided by their corresponding hyperbolic trigonometric denominators, yield zeta values (or $\log2$) of the correct ``degree".\medskip

\textsc{For Case~1:} when $n=a_1+a_2$ is odd, it suffices to show that the integrals 
\begin{align}
\int_0^\infty
\frac{\ell^{2k}\;\mathrm{d}\ell}
{\sinh(\ell)}
\end{align} 
are positive rational multiples of $\zeta(2k+1)$ for $k\geqslant 1$, and this follows from the following well-known formula (see \cite[Equation~3.523(1)]{MR2360010}):
\begin{align}
\label{eq:int-even-power-over-sinh}
\int_0^\infty
\frac{\ell^{2k}\;\mathrm{d}\ell}{\sinh(\ell)}
= 2(2k)!\left(1-\frac{1}{2^{2k+1}}\right)\zeta(2k+1). 
\end{align}

\textsc{For Cases~2 and 3:} when $n=a_1+a_2$ is even, it suffices to show that the integrals
\begin{align}
\int_0^\infty
\frac{\ell^{2k+1}\;\mathrm{d}\ell}{\sinh^2(\ell)}
,\quad 
\int_0^\infty\frac{\ell^{2k+1}\;\mathrm{d}\ell}{\cosh^2(\ell)}
\end{align}
are positive rational multiples of $\zeta(2k+1)$ for $k \geqslant 1$, and also that
\begin{align*}
\int_0^\infty
\frac{\ell\;\mathrm{d}\ell}{\cosh^2(\ell)}
\in\QQ_{>0}\cdot\log2.    
\end{align*}
This follows from the following well-known identities (see \cite[Equation~3.527(4),(3),(1)]{MR2360010}):
\begin{align*}
\int_0^\infty
\frac{\ell\;\mathrm{d}\ell}{\cosh^2(\ell)}
&= \log2,
\\
\int_0^\infty
\frac{\ell^{2k+1}\;\mathrm{d}\ell}{\cosh^2(\ell)}
&= \frac{(2k+1)!}{2^{2k}}(1-2^{-2k})\zeta(2k+1)
\,\quad \text{for }k\geqslant1,
\\
\int_0^\infty
\frac{\ell^{2k+1}\;\mathrm{d}\ell}{\sinh^2(\ell)} 
 &= \frac{(2k+1)!}{2^{2k}}\zeta(2k+1) 
\,\quad \text{for }k\geqslant1.
\end{align*}
\end{proof}

\subsection{Mirzakhani volumes of moduli spaces of general surfaces}

\begin{theorem}
\label{thm:generalmirzvolumes}
For $\Sigma\neq\DD_n, \AA_n$ or $\mathbb{A}_{a_1,a_2}$, let $a_{\max}:=\max_{1\leqslant i\leqslant l}\{a_i\}$. Then, the Mirzakhani volume $V_\Sigma(\vec{b})$ of $\mathcal{M}_\Sigma(\vec{b})$ is an element of
\begin{align*}
\QQ_{>0}[\pi^2, b_1^2,\cdots,b_m^2, \zeta(2j+1),\beta(2k)]_{j,k}
\end{align*}
with 
\begin{align*}
j,k\geqslant1;\quad
2j+1,2k \leqslant 6g-6+2(m+l)+ a_{\max}+1,    
\end{align*}
and where $\zeta(s)$ is the Riemann zeta function and $\beta(s)$ is the Dirichlet beta function. Moreover, $V_\Sigma(\vec{b})$ is a homogeneous polynomial of degree $\dim_{\mathbb{R}}\mathcal{M}_{\Sigma}(\vec{b})$, provided that $\pi$ is set to be degree $1$, each $b_i$ is degree~1, $\zeta(2j+1)$ is degree $2j+1$, and $\beta(2k)$ is degree $2k$. The degrees $\deg_{b_i}V_{\Sigma}(\vec{b})$ are at most $6g-6+2(m+l).$
\end{theorem}

\begin{proof}
We prove this by induction on the number of necks $l$. The $l=0$ base case is due to Mirzakhani \cite[Theorem~1.1]{mirz_simp}.\medskip

For the induction step: suppose that the statement holds for some $l\geqslant 0$, we prove it for $l+1$. Consider the neck $\nu_{l+1}$ of $\Sigma$, and let $\Sigma'$ be the component of $\Sigma\setminus\nu_{l+1}$ left after removing the $(l+1)$-th crown (which is topologically $\AA_{a_{l+1}}$). By the induction assumption,
\[
V_{\Sigma'}(b_1,\dotsc,b_m,\ell) \in \QQ_{>0}[\pi^2,b_1^2,\cdots,b_m^2,\ell^2, \zeta(2j+1),\beta(2k)]_{j,k},
\]
with $j,k\geqslant1$ and $2j+1,2k \leqslant 6g-6+2(m+1+l)+\max_{1\leqslant i\leqslant l}\{a_i\}+1$, and where $V_{\Sigma'}$ is of degree
\[
\dim_{\mathbb{R}}\mathcal{M}_{\Sigma'}(\vec{b},\ell) =6g-6+2(m+1)+3l+\sum_{k=1}^l a_k.
\]
Moreover, the degree $\deg_{\ell} V_{\Sigma'}(\vec{b},\ell)$ is at most $6g-6+2(m+1+l)$.

By \Cref{prop:vol-form-scc-new}:
\begin{align}
V_{\Sigma}(\vec{b}) = \int_0^\infty V_{\Sigma'}(b_1,\dotsc,b_m,\ell)\cdot  
\ell \cdot V_{\mathbb{A}_{a_{l+1}}}(\ell)\;\mathrm{d}\ell.
\end{align}

We split into two cases.\medskip

\noindent\textsc{Case~1:} $a_{l+1}$ is even. By \Cref{thm:fixedneckvol}, the expression
\[
\ell \cdot V_{\mathbb{A}_{a_{l+1}}}(\ell)\cdot  \sinh (\ell/2)
\]
is a homogeneous polynomial in $\QQ_{>0}[\pi^2,\ell^2]$ of degree $a_{l+1}$, and hence
\[
V_{\Sigma'}(b_1,\dotsc,b_m,\ell)\cdot  \ell \cdot V_{\mathbb{A}_{a_{l+1}}}(\ell) \cdot \sinh(\ell/2)
\]
is a homogeneous polynomial of degree $\dim_{\mathbb{R}}\mathcal{M}_{\Sigma'}(\vec{b},\ell) +a_{l+1}$, whose degree in $\ell$ is at most $6g-6+2(m+1+l)+a_{l+1}$. Further, by \Cref{eq:int-even-power-over-sinh} for $k\geqslant1$,
\begin{align}
\int_0^\infty
\frac{\ell^{2k}\;\mathrm{d}\ell}{\sinh(\ell/2)} \in \QQ_{>0} \cdot \zeta(2k+1), 
\end{align}
hence $V_{\Sigma}(\vec{b})$ is a homogeneous polynomial of degree
\[
\dim_{\mathbb{R}}\mathcal{M}_{\Sigma'}(\vec{b},\ell) +a_{l+1}+1 = 6g-6+2m+3(l+1)+\sum_{k=1}^{l+1} a_k = \dim_{\mathbb{R}}\mathcal{M}_{\Sigma}(\vec{b}),
\]
and moreover, 
\[
V_{\Sigma}(\vec{b}) \in \QQ_{>0}[\pi^2, b_1^2,\cdots,b_m^2, \zeta(2j+1),\zeta(a_{l+1}+1),\beta(2k)]_{j,k},
\]
with $j,k\geqslant1$ and $2j+1,2k \leqslant  6g-6+2(m+l+1)+\max_{1\leqslant i\leqslant l+1}\{a_i\}+1$.\\

\noindent\textsc{Case~2:} $a_{l+1}$ is odd. By \Cref{thm:fixedneckvol}, the expression
\[
\ell \cdot V_{\mathbb{A}_{a_{l+1}}}(\ell) \cdot \cosh (\ell/2)
\]
is a homogeneous polynomial in $\QQ_{>0}[\pi^2,\ell]$ which is \textit{odd} with respect to $\ell$ and has degree $a_{l+1}$. Recall \cite{wikipedia:dirichlet_beta} that the Dirichlet beta function satisfies
\[
\beta(s) 
:= \sum_{n=0}^\infty \frac{(-1)^n}{(2n+1)^s} 
= \frac{1}{\Gamma(s)}  \int_0^\infty\frac{x^{s-1}e^{-x}}{1+e^{-2x}}\mathrm{d}x,
\quad\text{where }\mathrm{Re}(s)>0.
\]    
Then,
\begin{equation}
\label{eq:dirichlet-beta}
\int_0^\infty\frac{\ell^{2k-1}}{\cosh(\ell/2)}\mathrm{d}\ell 
= 2^{2k+1}(2k-1)!\beta(2k)
\in \QQ_{>0}\cdot\beta(2k),
\end{equation}
for $k\geqslant 1$, and hence $V_{\Sigma}(\vec{b})$ is a homogeneous polynomial of degree
\[
\dim_{\mathbb{R}}\mathcal{M}_{\Sigma'}(\vec{b},\ell) +a_{l+1}+1 = \dim_{\mathbb{R}}\mathcal{M}_{\Sigma}(\vec{b}),
\]
and moreover, 
\[
V_{\Sigma}(\vec{b}) \in \QQ_{>0}[\pi^2, b_1^2,\cdots,b_m^2, \zeta(2j+1),\beta(2k),\beta(a_{l+1}+1)]{j,k},
\] 
with $j,k\geqslant1$ and $2j+1,2k \leqslant  6g-6+2(m+l+1)+\max_{1\leqslant i\leqslant l+1}\{a_i\}+1$. The proof is complete.
\end{proof}

\section{Volumes of moduli spaces of n-gons}
\label{subsec: ideal polygons}
In \cite[Proposition~6.2]{chekhov2024}, Chekhov obtains an expression for the Mirzakhani volume of $\mathcal{M}_{\mathbb{D}_n}$ as an integral over an $(n-3)$-dimensional simplex. Moreover, he explicitly carries out the volume calculation for the moduli spaces of $5$-gons and $6$-gons. We use \cite{MO402985} to show:\medskip 

\begin{theorem}
\label{thm:ngonvol}
The Mirzakhani volume of the moduli space of ideal $n$-gons, for $n\geqslant 3$, is
\begin{align*}
V_{\DD_n} = \left\{
\begin{array}{rl}
\frac{2}{n-2}\cdot\frac{(n-4)!!}{(n-3)!!}\pi^{n-4},
&\quad\text{if }n\text{ is even,}
\\[0.25em]
  \frac{1}{n-2}\cdot\frac{(n-4)!!}{(n-3)!!}\pi^{n-3},
&\quad\text{if }n\text{ is odd.}
\end{array}
\right.    
\end{align*}  
\end{theorem}

\begin{proof}
We take as starting point the integral expression for $V_{\mathbb{D}_n}$ given by Chekhov (implicitly) in \cite[Proposition~6.2]{chekhov2024}:

\begin{align*}
    V_{\mathbb{D}_n}
    =
    \int_{0<z_2<\cdots<z_{n-2}\leqslant1}
    \frac{\mathrm{d}z_2\wedge\ldots\wedge\mathrm{d}z_{n-2}}
    {(z_3-z_1)(z_4-z_2)\cdots(z_{n-2}-z_{n-4})(z_{n-1}-z_{n-3})},
\end{align*}
where $z_1=0$ and $z_{n-1}=1$. Take the change of coordinates given by $y_j=\frac{z_j}{z_{j+1}}$, then
\begin{align*}
\mathrm{d}y_2\wedge\cdots\wedge\mathrm{d}y_{n-2}
&=
(\tfrac{\mathrm{d}z_2}{z_3}-\tfrac{z_2\mathrm{d}z_3}{z_3^2})
\wedge\ldots\wedge
(\tfrac{\mathrm{d}z_{n-3}}{z_{n-2}}-\tfrac{z_{n-3}\mathrm{d}z_{n-2}}{z_{n-2}^2})
\wedge
\mathrm{d}z_{n-2}\\
&=
\frac{\mathrm{d}z_2\wedge\ldots\wedge\mathrm{d}z_{n-2}}
{z_3z_4\cdots z_{n-3} z_{n-2}}
\end{align*}
the above integral transforms to
\begin{align*}
V_{\mathbb{D}_n}
=\int_{(0,1)^{n-3}} \frac{\mathrm{d}y_2\wedge\cdots\wedge \mathrm{d}y_{n-2}}
{(1-\frac{z_1}{z_3})\cdots(1-\frac{z_{n-3}}{z_{n-1}})}
= \int_{(0,1)^{n-3}} \frac{\mathrm{d}y_2\wedge\cdots\wedge \mathrm{d}y_{n-2}}{\prod_{i=2}^{n-3}(1-y_iy_{i+1})}.
\end{align*}
Next, set $y_i = \frac{x_i-1}{x_i+1}$. Then $\mathrm{d}y_i = \frac{2}{(x_i+1)^2} \mathrm{d}x_i$, and
\begin{align*}
    1-y_i y_{i+1} = 1 - \frac{(x_i-1)(x_{i+1}-1)}{(x_i+1)(x_{i+1}+1)} = \frac{2(x_i+x_{i+1})}{(x_i+1)(x_{i+1}+1)}.
\end{align*}
Thus the volume integral transforms to
\begin{align}
\label{eq:vol-n-gon-int}
    V_{\mathbb{D}_n}
=2\int_{(1,\infty)^{n-3}} \frac{\mathrm{d}x_2\wedge\cdots\wedge \mathrm{d}x_{n-2}}{(1+x_2)(x_2+x_3)\dotsc(x_{n-3}+x_{n-2})(x_{n-2}+1)}.
\end{align}
To compute the integral in \Cref{eq:vol-n-gon-int}, we follow the strategy outlined in \cite{MO402985} and furnish the necessary details.

Consider the integral operator $T$ that acts on $f: [1,\infty) \to \RR$ by
\begin{align}
   T[f](y) = \int_1^{\infty} \frac{f(x)}{x+y} \mathrm{d}x,
\end{align}
provided that the integral exists for all $y\in [1,\infty)$. It follows from \Cref{eq:vol-n-gon-int} and the Fubini–Tonelli theorem that
\begin{align}
\label{eq:vol-n-gon-as-int-transform}
    V_{\mathbb{D}_n}
=2 \cdot T^{n-3}\left[\frac{1}{1+x}\right](1).
\end{align}
To clarify: to use Fubini-Tonelli, we assume that $V_{\DD_n}$ is finite. Chekhov asserts this in \cite{chekhov2024}, and we show it explicitly in \Cref{lem:fubini-tonelli}.
The key now is to find an integral transform that diagonalizes $T$. The Mehler-Fock transform \cite[Section~1.9]{MR2254107} is defined as 
\begin{align}
    \mathcal{MF}[f](\xi) = \int_1^\infty f(x)\cdot P_{-1/2+i\xi}(x)\,\mathrm{d}x,
\end{align}
where $P_{\nu}$ is the Legendre function of the first kind \cite[Section~7.3]{MR350075} of degree $\nu \in \CC$. Sufficient conditions for the existence of the Mehler-Fock transform are \begin{itemize}
    \item $f\in L^{\text{loc}}_1(1,\infty)$, and 
    \item $f(x) = O(x^{\alpha}),\, x\to \infty$ for some $\alpha<-1/2$ \cite[Theorem~1.9.51]{MR2254107}. 
\end{itemize}
Both conditions are satisfied for $f(x) = \frac{1}{1+x}$. We invoke the following important identity due to Mehler \cite[Page~193]{MR1510098}:
\begin{align}
\label{eq: Mehler formula}
    \int_1^\infty \frac{1}{x+y} \cdot P_{-1/2+i\xi}(x) \,\mathrm{d}x = \frac{\pi}{\cosh \pi \xi}\cdot P_{-1/2+i\xi}(y),
    \text{ for all $y\geqslant1$.}
\end{align}
 Then, if the Fubini–Tonelli theorem applies,
\begin{equation}
\begin{split}
    \mathcal{MF}[Tf](\xi) &=  \int_1^\infty\left[\int_1^{\infty} \frac{f(x)}{x+y} \mathrm{d}x\right] \cdot P_{-1/2+i\xi}(y)\,\mathrm{d}y 
    \\& 
    = \int_1^\infty f(x) \left[\int_1^\infty \frac{1}{x+y}\cdot P_{-1/2+i\xi}(y) \,\mathrm{d}y\right ] \mathrm{d}x
    \\&
    =\frac{\pi}{\cosh \pi\xi}\cdot \mathcal{MF}[f](\xi),
\end{split}
\end{equation}
and more generally
\begin{align*}
    \mathcal{MF}[T^kf](\xi) = \frac{\pi^k}{\cosh^k \pi\xi} \mathcal{MF}[f](\xi),
    \text{ for $k\geqslant 1$.}
\end{align*}
In particular, \Cref{eq: Mehler formula} and \cite[Equation~7.3.13]{MR350075} implies that 
\begin{align*}
    \mathcal{MF}\left[\frac{1}{1+x}\right](\xi) = \frac{\pi}{\cosh \pi \xi} \cdot P_{-1/2+i\xi}(1) = \frac{\pi}{\cosh \pi \xi}.
\end{align*}
By \Cref{prop:fubini-tonelli}, Fubini--Tonelli applies and hence
\begin{align}
\label{eq:MF-of-iterates-of-f}    \mathcal{MF}\left[T^{k}\left[\frac{1}{1+x}\right]\right](\xi) = \frac{\pi^{k+1}}{\cosh^{k+1}\pi\xi},
\text{ for $k\geqslant 1$.}\end{align}
 The inverse Mehler-Fock transform is given by \cite[Equation~1.9.11]{MR2254107}
\begin{align}
\label{eq:MF-inverse}
    \mathcal{MF}^{-1}[F](x) = \int_0^\infty \xi \tanh\pi\xi \cdot F(\xi)P_{-1/2+i\xi}(x)\,\mathrm{d}\xi.
\end{align}
Combining \Cref{eq:vol-n-gon-as-int-transform}, \Cref{eq:MF-of-iterates-of-f} and \Cref{eq:MF-inverse}, we obtain  
\begin{equation}
\begin{split}
V_{\DD_{n}} &= 2 \cdot T^{n-3}\left[\frac{1}{1+x}\right](1) = 2\cdot \mathcal{MF}^{-1}\left[\frac{\pi^{n-2}}{\cosh^{n-2}\pi\xi}\right](1)
\\& 
= 2 \cdot \int_0^\infty \xi \tanh\pi\xi  \cdot \frac{\pi^{n-2}}{\cosh^{n-2}\pi\xi} \cdot P_{-1/2+i\xi}(1)\,\mathrm{d}\xi 
\\&
= 2\pi^{n-2}\cdot \int_0^\infty  \frac{\xi \sinh\pi\xi}{\cosh^{n-1}\pi\xi}\mathrm{d}\xi.
\end{split}
\end{equation}
We proceed via integration by parts. Let 
\begin{align*}
u = \xi,\quad v = \int \frac{ \sinh\pi\xi}{\cosh^{n-1}\pi\xi}\mathrm{d}\xi = -\frac{1}{(n-2)\pi}\cosh^{2-n}(\pi\xi).
\end{align*}
Then
\begin{align*}
V_{\DD_{n}} &=   2\pi^{n-2} \int_0^\infty u\cdot v'\, \mathrm{d}\xi = 2\pi^{n-2} \left( u v\Big|_{0}^\infty -\int_0^\infty u'\cdot v\,\mathrm{d}\xi\right)
\\& = 2\pi^{n-2} \left(- \frac{\xi}{(n-2)\pi}\cosh^{2-n}(\pi\xi)\Big|_{0}^\infty +\int_0^\infty  \frac{1}{(n-2)\pi}\cosh^{2-n}(\pi\xi)\; \mathrm{d}\xi \right)
\\& = \frac{2\pi^{n-3}}{n-2}\int_0^\infty  \cosh^{2-n}(\pi\xi)\;\mathrm{d}\xi,
\end{align*}
where in the last line we used that $n>2$. Invoking \cite[Equation~3.512(2)]{MR2360010} yields
\begin{align}
V_{\DD_{n}}& = \frac{2\pi^{n-3}}{n-2} \cdot \frac{1}{2\pi} B\left(\frac{1}{2},\frac{n}{2}-1\right) = \frac{\pi^{n-4}}{n-2} \cdot \frac{\Gamma\left(\frac{1}{2}\right)\Gamma\left(\frac{n-2}{2}\right)}{\Gamma\left(\frac{n-1}{2}\right)},
\end{align}
where $B$ is the beta function.\medskip

If $n=2k$ is even, invoking \cite[Equation~1.2(15)]{batemanHTF1} yields  
\begin{equation*}
\begin{split}
V_{\DD_{n}} &= \frac{\pi^{n-4}}{n-2} \cdot\frac{\Gamma\left(\frac{1}{2}\right)\Gamma\left(\frac{n-2}{2}\right)}{\Gamma\left(\frac{n-1}{2}\right)} =\frac{\pi^{2k-4}}{2k-2} \cdot\frac{\Gamma\left(\frac{1}{2}\right)\Gamma\left(k-1\right)}{\Gamma\left(k-\frac{1}{2}\right)} 
\\& = \frac{\pi^{2k-4}}{2k-2} \cdot\frac{\sqrt{\pi}(k-2)!}{\frac{(2k-3)!!}{2^{k-1}}\sqrt{\pi}} = \frac{\pi^{2k-4}}{2k-2} \cdot\frac{2(2k-4)!!}{(2k-3)!!}
\\& = \frac{2}{n-2} \cdot \frac{(n-4)!!}{(n-3)!!}\pi^{n-4}.
\end{split}
\end{equation*}
If $n=2k+1$ is odd, invoking \cite[Equation~1.2(15)]{batemanHTF1} yields
\begin{equation*}
\begin{split}
V_{\DD_{n}} &= \frac{\pi^{n-4}}{n-2} \cdot\frac{\Gamma\left(\frac{1}{2}\right)\Gamma\left(\frac{n-2}{2}\right)}{\Gamma\left(\frac{n-1}{2}\right)} =\frac{\pi^{2k-3}}{2k-1} \cdot\frac{\Gamma\left(\frac{1}{2}\right)\Gamma\left(k-\frac{1}{2}\right)}{\Gamma\left(k\right)} 
\\& = \frac{\pi^{2k-3}}{2k-1} \cdot\frac{\sqrt{\pi}\cdot\frac{(2k-3)!!}{2^{k-1}}\sqrt{\pi}}{(k-1)!} = \frac{\pi^{2k-2}}{2k-1} \cdot\frac{(2k-3)!!}{(2k-2)!!}
\\& = \frac{1}{n-2} \cdot \frac{(n-4)!!}{(n-3)!!}\pi^{n-3}.
\end{split}
\end{equation*}

\end{proof}

\subsection{Generating function}
We next derive a generating function for the volumes of moduli spaces of $n$-gons. Since these moduli spaces only start making sense for $n\geqslant3$, we set the generating function to take the form
\[
P_{\mathbb{D}_{\bullet}}(x):=\sum_{n=3}^\infty V_{\mathbb{D}_n}\,x^{n-3}.
\]

\begin{proposition}[generating function for volumes of moduli spaces of $n$-gons]
\label{prop:genfunction}
The generating function $P_{\mathbb{D}_{\bullet}}(x)$ for the Mirzakhani volumes of moduli spaces of $n$-gons is
\begin{align}
P_{\mathbb{D}_{\bullet}}(x)
= 
\frac{\arcsin{\pi x}}{\pi x}
+ 
x\Bigl(\frac{\arcsin{\pi x}}{\pi x}\Bigr)^2. 
\end{align}
\end{proposition}

\begin{proof}
We first note that $\frac{\arcsin{\pi x}}{\pi x}$ is an even function and $x\left(\frac{\arcsin{\pi x}}{\pi x}\right)^2$ is an odd function, and hence they respectively correspond to the odd $n$ volume terms and even $n$ volume terms in the generating function. Now, note that the Taylor series expansion of $\arcsin{x}$ at $x=0$ is \cite[Equation~1.641(1)]{MR2360010}:
\[
\arcsin{x} = \sum_{j=0}^{\infty} \frac{(2j-1)!!}{(2j)!!}\frac{x^{2j+1}}{2j+1}
\quad\Rightarrow\quad
 \frac{\arcsin{\pi x}}{\pi x} 
 = \sum_{j=0}^{\infty} \frac{(2j-1)!! \pi^{2n}}{(2j)!!}\frac{x^{2j}}{2j+1}.
\]
Which yields the desired volumes for $n-3=2j$. Similarly \cite[Equation~1.645(2)]{MR2360010}:
\[
(\arcsin{x})^2 
= \sum_{j=0}^{\infty} 
\frac{2\cdot(2j)!!}{(2j+2)\cdot(2j+1)!!} x^{2j+2},
\]
and hence
\[
x\left(\frac{\arcsin{\pi x}}{\pi x}\right)^2 
= 
\sum_{j=0}^{\infty} 
\frac{2\cdot(2j)!!\;\pi^{2j}}{(2j+2)\cdot(2j+1)!!} x^{2j+1}. 
\]
which are precisely the desired volumes for $n-3=2j+1$.
\end{proof}

\section{Further exploration}

\subsection{Closed formulae}

In the course of our investigations, we have generated Mirzakhani volumes such as the following for the moduli space of $\mathbb{A}_{a_1,a_2}$:

\begin{center}
\begin{tabular}{ |l|l| } 
\hline
    $\Sigma$ & $V_\Sigma$\\ \hline
    $\mathbb{A}_{1,1}$ & $\log 2$  \\ \hline
    $\mathbb{A}_{1,2}$  & $\frac{7\zeta(3)}{4}$  \\ \hline
    $\mathbb{A}_{1,3}$ & $\frac{\pi ^2 \log 2}{2} +\frac{9\zeta (3)}{4} $    \\ \hline
    $\mathbb{A}_{1,4}$ & $\frac{7\pi ^2 \zeta (3)}{6}+\frac{31 \zeta (5)}{8}$   \\ \hline
    $\mathbb{A}_{1,5}$ & $
    \frac{3\pi^4\log 2}{8}+\frac{15\pi^2\zeta(3)}{8}+\frac{75\zeta(5)}{16}$  \\ \hline
    $\mathbb{A}_{1,6}$ & $\frac{14 \pi ^4 \zeta (3)}{15}+\frac{31 \pi ^2 \zeta (5)}{8}+\frac{381 \zeta (7)}{64}$   \\ \hline
    $\mathbb{A}_{1,7}$ & $\frac{5\pi^6\log 2}{16}+\frac{259\pi^4\zeta(3)}{160}+\frac{175\pi^2\zeta(5)}{32}+\frac{441\zeta(7)}{64}$  \\ \hline
    $\mathbb{A}_{1,8}$ & $\frac{4\pi^6\zeta(3)}{5}+\frac{217\pi^4\zeta(5)}{60}+\frac{127\pi^2\zeta(7)}{16}+\frac{511\zeta(9)}{64}$  \\ \hline
    $\mathbb{A}_{1,9}$ & $\frac{35\pi^8\log2}{128}+\frac{3229\pi^6\zeta(3)}{2240}+\frac{705\pi^4\zeta(5)}{128}+\frac{1323\pi^2\zeta(7)}{128}+\frac{2295\zeta(9)}{256}$  \\ \hline
    $\mathbb{A}_{1,10}$ & $\frac{32 \pi^8 \zeta(3)}{45} + \frac{1271 \pi^6 \zeta(5)}{378} + \frac{1651 \pi^4 \zeta(7)}{192} + \frac{2555 \pi^2 \zeta(9)}{192} + \frac{10235 \zeta(11)}{1024}$  \\ \hline
    $\mathbb{A}_{1,11}$ & $\frac{63 \pi^{10} \log(2)}{256} + \frac{117469 \pi^8 \zeta(3)}{89600} + \frac{86405 \pi^6 \zeta(5)}{16128} + \frac{30723 \pi^4 \zeta(7)}{2560} + \frac{8415 \pi^2 \zeta(9)}{512} + \frac{11253 \zeta(11)}{1024}$    \\ \hline
    $\mathbb{A}_{1,12}$ & $\frac{64 \pi^{10} \zeta(3)}{99} + \frac{14849 \pi^8 \zeta(5)}{4725} + \frac{17653 \pi^6 \zeta(7)}{2016} + \frac{15841 \pi^4 \zeta(9)}{960} + \frac{10235 \pi^2 \zeta(11)}{512} + \frac{24573 \zeta(13)}{2048}$    \\ \hline
    $\mathbb{A}_{2,2}$ & $6\zeta(3)$   \\ \hline
    $\mathbb{A}_{2,3}$ & $\frac{7\pi ^2 \zeta (3)}{8}+\frac{93 \zeta (5)}{8}$    \\ \hline
    $\mathbb{A}_{2,4}$ & $4\pi^2\zeta(3)+20\zeta(5)$   \\ \hline
    $\mathbb{A}_{2,5}$ & $\frac{21 \pi ^4 \zeta (3)}{32}+\frac{155 \pi ^2 \zeta (5)}{16}+\frac{1905 \zeta (7)}{64}$   \\ \hline
    $\mathbb{A}_{3,3}$ & $
    \frac{\pi^4\log 2}{4}+\frac{9\pi^2\zeta(3)}{4}+\frac{225\zeta(5)}{8}$  \\ \hline
    $\mathbb{A}_{3,4}$ & $\frac{7 \pi ^4 \zeta (3)}{12}+\frac{155 \pi ^2 \zeta (5)}{16}+\frac{1905 \zeta (7)}{32}$   \\ \hline
\end{tabular}
\end{center}

\begin{remark}
    
Based on these and other computed exact volumes, we have observed that the constants in these formulae seem to satisfy the following (unverified) patterns:

\begin{enumerate}
\item
The leading term of the volume of the moduli space of $(1,2k+1)$-annuli, when regarded as a polynomial in $\pi^2$, is as follows:
\[
\frac{(2k-1)!!}{(2k)!!}\pi^{2k}\log2.
\]
The leading term of the volume of the moduli space of $(1,2k)$-annuli is as follows:
\[
\frac{7}{4}\frac{(2k-2)!!}{(2k-1)!!}\pi^{2k-2}\zeta(3).
\]
\item 
The last term of the volume of the moduli space of $(1,k)$-annuli (for $k\geqslant 2$) is
    \begin{align*}
    k(1-2^{1-k})\zeta(2\lfloor \tfrac{k}{2}\rfloor+1),\quad \text{ when }k\text{ is odd}\\
    k(1-2^{-1-k})\zeta(2\lfloor \tfrac{k}{2}\rfloor+1),\quad \text{ when }k\text{ is even.}
    \end{align*}

\item The second last term of the volume of the moduli space of $(1,k)$-annuli (for $k\geqslant 4$) is
    \begin{align*}
    \tfrac{k(k-2)}{3!}(1-2^{3-k})\zeta(2\lfloor \tfrac{k}{2}\rfloor-1)\pi^2,\quad \text{ when }k\text{ is odd}\\
    \tfrac{k(k-2)}{3!}(1-2^{1-k})\zeta(2\lfloor \tfrac{k}{2}\rfloor-1)\pi^2,\quad \text{ when }k\text{ is even.}
    \end{align*}
\item 
The third last term of the volume of the moduli space of $(1,k)$-annuli (for $k\geqslant 6$) is
    \begin{align*}
    \tfrac{2k(k-4)(5k+2)}{6!}(1-2^{5-k})\zeta(2\lfloor \tfrac{k}{2}\rfloor-3)\pi^4,\quad \text{ when }k\text{ is odd}\\
    \tfrac{2k(k-4)(5k+2)}{6!}(1-2^{3-k})\zeta(2\lfloor \tfrac{k}{2}\rfloor-3)\pi^4,\quad \text{ when }k\text{ is even.}
    \end{align*}
\item 
The fourth last term of the volume of the moduli space of $(1,k)$-annuli (for $k\geqslant 8$) is
    \begin{align*}
    \tfrac{8k(k-6)(35k^2+42k+16)}{9!}(1-2^{7-k})\zeta(2\lfloor \tfrac{k}{2}\rfloor-5)\pi^6,\quad \text{ when }k\text{ is odd}\\
    \tfrac{8k(k-6)(35k^2+42k+16)}{9!}(1-2^{5-k})\zeta(2\lfloor \tfrac{k}{2}\rfloor-5)\pi^6,\quad \text{ when }k\text{ is even.}
    \end{align*}
\item 
The fifth last term of the volume of the moduli space of $(1,k)$-annuli (for $k\geqslant 10$) is
    \begin{align*}
    \tfrac{88k(k-8)(5k+4)(35k^2+56k+36)}{12!}(1-2^{9-k})\zeta(2\lfloor \tfrac{k}{2}\rfloor-7)\pi^8,\quad \text{ when }k\text{ is odd}\\
    \tfrac{88k(k-8)(5k+4)(35k^2+56k+36)}{12!}(1-2^{7-k})\zeta(2\lfloor \tfrac{k}{2}\rfloor-7)\pi^8,\quad \text{ when }k\text{ is even.}
    \end{align*}
\item 
The sixth last term of the volume of the moduli space of $(1,k)$-annuli (for $k\geqslant 12$) is
    \begin{align*}
    \tfrac{3640k(k-10)(385k^4+1540k^3+2684k^2+2288k+768)}{15!}(1-2^{11-k})\zeta(2\lfloor \tfrac{k}{2}\rfloor-9)\pi^{10},\; \text{ when }k\text{ is odd}\\
    \tfrac{3640k(k-10)(385k^4+1540k^3+2684k^2+2288k+768)}{15!}(1-2^{9-k})\zeta(2\lfloor \tfrac{k}{2}\rfloor-9)\pi^{10},\; \text{ when }k\text{ is even.}
    \end{align*}
\end{enumerate}
\end{remark}
\medskip

Recall \cite[Section~1.3]{MR2868112} that Stirling numbers of the first kind are defined via
\begin{align}
\label{eq:stirling-numbers}
    x(x-1)\dotsc(x-n+1) = \sum_{k=0}^n s(n,k)x^k.
\end{align}
We prove:
\begin{proposition}
\label{prop:coeff-vol-bi-annuli}
For $a_1=1, a_2 = 2k$, the coefficients $\rho_{2i}(a_1,a_2)$ in \Cref{thm:vol-biannuli} satisfy
\begin{align}
\label{eq:coeff-vol-bi-annuli}
\rho_{2i}(1,2k)  = \frac{2^{2i}(1-\frac{1}{2^{2k+1-2i}})(2k-2i)!}{(2k-1)!} \cdot\sum_{r=0}^{2k-2i}(-1)^{k-i+r} s(k,r)s(k,2k-2i-r). 
\end{align}
\end{proposition}

\begin{proof}
By \Cref{eq:int-vol-bi-annulus} and \Cref{thm:fixedneckvol},
\begin{equation}
\begin{split}
\label{eq:vol-1-2k-annulus-int}
V_{\mathbb{A}_{1,2k}} 
&= \int_0^\infty 
\ell
\cdot V_{\mathbb{A}_{1}}(\ell)
\cdot V_{\mathbb{A}_{2k}}(\ell) 
\,\mathrm{d}\ell \\&=  \int_0^\infty \ell\cdot \frac{1}{2\cosh \ell/2}\cdot  \frac{\ell}{\sinh{\ell/2}} \cdot 
\frac{\prod_{j=1}^{k-1} (\ell^2+(2j)^2\pi^2)}{2(2k-1)!}\,\mathrm{d}\ell
\\& 
= \int_0^\infty \frac{\ell^2 \prod_{j=1}^{k-1} (\ell^2+(2j)^2\pi^2)}{2(2k-1)!\sinh \ell } \,\mathrm{d}\ell.
\end{split}
\end{equation}
Next, we expand the numerator of the integrand in \Cref{eq:vol-1-2k-annulus-int}. Note that 
\begin{align}
\label{eq:sum-of-two-squares}
\ell^2+(2j)^2\pi^2 =\left(\ell+(2j)\pi\sqrt{-1}\right)\left(\ell-(2j)\pi\sqrt{-1} \right),
\end{align}
and also that by \Cref{eq:stirling-numbers},
\begin{align}
\label{eq:prod-poly-plus}
\ell\cdot\prod_{j=1}^{k-1} \left(\ell+(2j)\pi\sqrt{-1}\right) = \sum_{j=0}^{k}(-1)^{k-j} s(k,j)\left(2\pi\sqrt{-1}\right)^{k-j}\ell^j, 
\end{align}
\begin{align}
\label{eq:prod-poly-minus}
\ell\cdot\prod_{j=1}^{k-1} \left(\ell-(2j)\pi\sqrt{-1}\right) = \sum_{j=0}^{k} s(k,j)\left(2\pi\sqrt{-1}\right)^{k-j}\ell^j. 
\end{align}
Since the numerator of the integrand in \Cref{eq:vol-1-2k-annulus-int} is an even polynomial, together with \Cref{eq:sum-of-two-squares,eq:prod-poly-plus,eq:prod-poly-minus}, we obtain: 
\[
\begin{aligned}
&\ell^2\cdot\prod_{j=1}^{k-1} (\ell^2+(2j)^2\pi^2) \\
\smash{\mathllap{=}}\, &\sum_{j=0}^{k} \left( \sum_{r=0}^{2j}(-1)^{k-r} s(k,r)s(k,2j-r)\left(2\pi\sqrt{-1}\right)^{2k-2j} \right)  \ell^{2j}\\
\smash{\mathllap{=}}\, & \sum_{j=0}^{k} \left( \sum_{r=0}^{2j}(-1)^{j+r} s(k,r)s(k,2j-r)\left(2\pi\right)^{2k-2j} \right)  \ell^{2j}.
\end{aligned}
\]
Recall from \Cref{eq:int-even-power-over-sinh} that for $j\geqslant1$,
\begin{align}
\int_0^\infty
\frac{\ell^{2j}\;\mathrm{d}\ell}{\sinh(\ell)}
= 2(2j)!\left(1-\frac{1}{2^{2j+1}}\right)\zeta(2j+1). 
\end{align}
Hence we compute
\begin{equation}
\label{eq:vol-1-2k-annulus-sum}
\begin{split}
V_{\mathbb{A}_{1,2k}} 
&= \int_0^\infty \frac{\ell^2 \prod_{j=1}^{k-1} (\ell^2+(2j)^2\pi^2)}{2(2k-1)!\sinh \ell } \,\mathrm{d}\ell
\\&
= \frac{1}{2(2k-1)!} \sum_{j=0}^{k} \left( \sum_{r=0}^{2j}(-1)^{j+r} s(k,r)s(k,2j-r)\left(2\pi\right)^{2k-2j} \right)  
\\&
\quad\cdot \left(2(2j)!\left(1-\frac{1}{2^{2j+1}}\right)\zeta(2j+1)\right) 
\\&
=\sum_{j=0}^{k}\left(  \frac{2^{2k-2j}\left(1-\frac{1}{2^{2j+1}}\right)(2j)!}{(2k-1)!} \sum_{r=0}^{2j}(-1)^{j+r} s(k,r)s(k,2j-r)\right)
\\&
\quad\cdot\pi^{2k-2j}\zeta(2j+1).
\end{split}
\end{equation}
Note that if $j=0$, then $s(k,2j)=0$. Then the proposition follows by letting $2i=2k-2j$ in \Cref{eq:vol-1-2k-annulus-sum}.

\end{proof}

\begin{remark} 
It is not immediate from \Cref{prop:coeff-vol-bi-annuli} alone that the coefficients $\rho_{2i}(1,2k)$ are positive. Indeed, not all terms in the summation in \Cref{eq:coeff-vol-bi-annuli} are positive, e.g. for $i=1,k=3$ we have
\begin{equation*}
\begin{split}
\rho_{2}(1,6) &= \frac{2^2(1-\frac{1}{2^5})4!}{5!}\cdot \left((-1)^{3}\cdot2\cdot1+(-1)^{4}\cdot(-3)\cdot(-3)+(-1)^{5}\cdot1\cdot2\right) \\&
=\frac{31}{40} (-2+9-2) = \frac{31}{8}.     
\end{split}    
\end{equation*}
\end{remark}

\begin{remark}
A double summation formula \cite[Equation~5.7a]{MR460128} for Stirling numbers of the first kind involving only products and quotients of factorials and powers is known (but no such formula with single summation is known!). It yields a formula for $\rho_{2i}(1,2k)$ involving \textit{five} ``nested'' summations.

Further, similar (but more cumbersome) formulas can be obtained for the volumes $V_{\AA_{2a_1,2a_2}}$. We did not pursue a similar description for general volumes $V_{\AA_{a_1,a_2}}$. 
\end{remark}

\newpage
\appendix

\newpage

\section{Volume computations via lambda-lengths}

\subsection{Volume of the moduli space of 2-crowns via lambda-lengths}
\label{app:2-crowns}

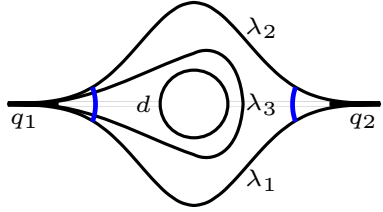
\begin{figure}[H]
    \centering
    \begin{tikzpicture}[scale=1.5]

    \draw[thick, smooth cycle] 
        plot[domain=-1:1] (\x, {0.6*exp(-5*\x*\x)}) -- 
       plot[domain=1:-1] (\x, {-0.6*exp(-5*\x*\x)}) -- cycle;

    \draw[very thick, white] (-0.8,0) -- (0.8,0);

    \filldraw[fill=white, thick] (0,0) circle (0.2cm);

     New lasso (black, around the eyeball, from and to the left corner)
    \draw[thick] 
        (-1,0) .. controls (-0.6,0) and (-0.3,0.2) .. (0,0.3)
              .. controls (0.1,0.35) and (0.27,0.3) .. (0.29,0)
            .. controls (0.27,-0.3) and (0.1,-0.35) .. (0,-0.3)
               .. controls (-0.3,-0.2) and (-0.6,0) .. (-1,0);

\draw (0.4,0.45) node {\tiny $\lambda_2$};

\draw (0.4,-0.45) node {\tiny $\lambda_1$};

\draw (0.4,0) node {\tiny $\lambda_3$};

\draw (-0.3,0) node {\tiny $d$};

\draw (-1,-0.1) node {\tiny $q_1$};

\draw (1,-0.1) node {\tiny $q_2$};

\draw[very thick, blue] (0.6,0.1) arc (160:200:0.3);

\draw[very thick, blue] (-0.6,0.1) arc (20:-20:0.3);

\end{tikzpicture}
    \caption{Crown with two tines, triangulated by an arc. In \textcolor{blue}{blue}: unit horocycles.}
    \label{fig:2-tines-lambda-lengths}
\end{figure}

We give an alternative proof of \Cref{thm:fixedneckvol} in the case $n=2$:
\[
V_{\mathbb{A}_{2}}(d) = \frac{d/2}{\sinh{d/2}},
\]
using $\lambda$-lengths. Consider the triangulation of the 2-crown by a single arc with endpoints at the boundary puncture $q_1$, and decorate it with unit horocycles as in \figref{fig:2-tines-lambda-lengths}. Let $\lambda_1, \lambda_2,  \lambda_3$ be the lambda lengths of the resulting truncated arcs as in \figref{fig:2-tines-lambda-lengths}. 
We express the horocycle lengths in terms of the lambda lengths (see \cite[Lemma~4.9]{pennerbook}, page 36, and \cite[Note~3.4]{huang_thesis} for the formulas we use):
\begin{equation}
\label{eq:2-tine-crown-hor-cond-1}
1 = \frac{2\cosh{d/2}}{\lambda_3} + \frac{\lambda_1}{\lambda_2\lambda_3}+\frac{\lambda_2}{\lambda_1\lambda_3},
\end{equation}
\begin{equation}
\label{eq:2-tine-crown-hor-cond-2}
1 = \frac{\lambda_3}{\lambda_1 \lambda_2}.
\end{equation}
Multiplying \Cref{eq:2-tine-crown-hor-cond-1} by $\lambda_1\lambda_3/\lambda_2$ and using that $\lambda_3 = \lambda_1 \lambda_2$ from \Cref{eq:2-tine-crown-hor-cond-2}, we get:
\begin{equation}
\label{eq:2-tine-crown-hor-cond-3}
\frac{2\lambda_1\cosh{d/2}}{\lambda_2} + \frac{\lambda_1^2}{\lambda_2^2}+1 = \frac{\lambda_1\lambda_3}{\lambda_2} = \lambda_1^2.
\end{equation}
Next, following \cite[Page~9]{chekhov2024}, the decoration-independent coordinate $x_1$ is defined as 
\[
x_1 = \frac{\lambda_1}{\lambda_2},
\]
thus we can rewrite \Cref{eq:2-tine-crown-hor-cond-3} as
\begin{equation}
\label{eq:2-tine-crown-hor-cond-4}
2x_1 \cosh d/2 + x_1^2+1 =  \lambda_1^2.
\end{equation} 
From \Cref{eq:2-tine-crown-hor-cond-4}, Chekhov's functional can be expressed as follows:
\begin{align*}
e^{-S} = \frac{1}{\lambda_1\lambda_2} = \frac{x_1}{\lambda_1^2} = \frac{x_1}{x_1^2+2x_1 \cosh d/2+1}.
\end{align*}
By \cite[Lemma~2.2]{chekhov2024}, the volume form on the moduli space of 2-crowns equals 
\begin{align*}
\Omega^{\mathrm{WP}}_{\AA_2}(d) = \mathrm{d} \log x_1 =\frac{\mathrm{d}x_1}{x_1}
\end{align*}
Putting it together, we have
\begin{equation}
\begin{split}
V_{\mathbb{A}_{2}}(d) &= \int e^{-S} \Omega^{\mathrm{WP}}_{\AA_2}(d) =  \int_0^\infty \frac{x_1}{x_1^2+2x_1 \cosh d/2+1} \cdot \frac{\mathrm{d}x_1}{x_1} \\& 
= \int_0^\infty \frac{\mathrm{d}x_1}{(x_1+\cosh d/2)^2-\sinh^2d/2}.
\end{split}
\end{equation}
Let $u = x_1+\cosh d/2; a= \sinh d/2$, then 
\begin{equation}
\begin{split}
V_{\mathbb{A}_{2}}(d) &=     \int_{\sqrt{a^2+1}}^\infty \frac{\mathrm{d}u}{u^2-a^2} = \int_{\sqrt{a^2+1}}^\infty\frac{1}{2a}\Bigl( \frac{1}{u-a}-\frac{1}{u+a} \Bigr) \mathrm{d}u \\& = \frac{1}{2a}\log\Bigl(\frac{u-a}{u+a}\Bigr)\Big|_{\sqrt{a^2+1}}^\infty
 = \frac{1}{2a}\log\Biggl(\frac{\sqrt{a^2+1}+a}{\sqrt{a^2+1}-a}\Biggr) \\& = \frac{1}{2\sinh d/2}\log \frac{e^{d/2}}{e^{-d/2}} =\frac{d/2}{\sinh d/2}.
\end{split}    
\end{equation}
\medskip

\section{Fubini-Tonelli, Mehler-Fock and n-gons}
For a fixed integer $n\geqslant4$ and arbitrary $y\geqslant 1$, define the function
\begin{align}
V_{\DD_n}(y) = 2\int_{(1,\infty)^{n-3}} \frac{\mathrm{d}x_2\wedge\cdots\wedge \mathrm{d}x_{n-2}}{(1+x_2)(x_2+x_3)\dotsc(x_{n-3}+x_{n-2})(x_{n-2}+y)}. 
\end{align}
Note that in particular $V_{\DD_{n}}(1)=V_{\DD_n}$. We prove:
\begin{proposition}
\label{prop:fubini-tonelli}
For every $y\geqslant 1,$ the integral $V_{\DD_n}(y)$ is finite. Moreover, for every $\xi\geqslant1$,
\begin{align}
    \int_1^\infty |V_{\DD_n}(y)\cdot P_{-1/2+i\xi}(y)|\,\mathrm{d}y <\infty.
\end{align}
\end{proposition}
We start with the following lemma:
\begin{lemma}
\label{lem:fubini-tonelli}
There exist constants $C_n,D_n\geqslant1$ such that the following holds:
\begin{equation}
\begin{split}
    V_{\DD_{n}}(y) \leqslant C_n,& \quad1\leqslant y\leqslant 2,
    \\
    V_{\DD_{n}}(y) \leqslant D_n \frac{\sum_{k=0}^{n-3}\log^k(y)}{y},& \quad 2\leqslant y < \infty.
\end{split}
\end{equation}
\end{lemma}
\begin{proof}
We prove this by induction on $n$. For $n=4$ base case:
\begin{equation}
\begin{split}
V_{\DD_4}(y) & =  2\int_1^\infty\frac{\mathrm{d}x}{(1+x)(x+y)} \\&= \frac{2}{y-1}\int_1^\infty \left(\frac{1}{1+x}-\frac{1}{x+y}\right)\mathrm{d}x \\&= \frac{2\log\frac{y+1}2{}}{y-1} \leqslant \frac{2\log y}{y-1},
\end{split}
\end{equation}
where in the last line we used that $\frac{y+1}{2}\leqslant y$ for $y\geqslant1$. Note that if $y\geqslant 2$, then $\frac{1}{y-1}\leqslant \frac{2}{y}$ and thus $\frac{2\log y}{y-1}\leqslant \frac{4\log y}{y}$. It is straightforward to check that $V_{\DD_4}(1)=1$ and that $V_{\DD_4}(y)$ is decreasing for $y\geqslant1$. Hence we can let $C_4=1,D_4=4$.

For the induction step: suppose that the statement holds for some $n$, we prove it for $n+1$. For a fixed $z\geqslant1$, consider the integral
\begin{equation}
\begin{split}
\int_1^\infty \frac{V_{\DD_n}(y)}{y+z} \mathrm{d}y 
= \int_1^{2} \frac{V_{\DD_n}(y)}{y+z} \mathrm{d}y+\int_{2}^z \frac{V_{\DD_n}(y)}{y+z} \mathrm{d}y+\int_z^\infty \frac{V_{\DD_n}(y)}{y+z} \mathrm{d}y.
\end{split}
\end{equation}
Since $\frac{1}{y+z}\leqslant \frac{1}{y}, \frac{1}{y+z}\leqslant \frac{1}{z}$, by the induction assumption we have 
\begin{equation}
\label{eq:int-3-terms}
\begin{split}
\int_1^\infty \frac{V_{\DD_n}(y)}{y+z} \mathrm{d}y  & 
\leqslant \int_1^{2} \frac{V_{\DD_n}(y)}{z} \mathrm{d}y+\int_{2}^z \frac{V_{\DD_n}(y)}{z} \mathrm{d}y+\int_z^\infty \frac{V_{\DD_n}(y)}{y} \mathrm{d}y
\\&
\leqslant \frac{C_n}{z} + \frac{1}{z}\int_{2}^z D_n \frac{\sum_{k=0}^{n-3}\log^k(y)}{y}\mathrm{d}y+ \int_z^\infty D_n\frac{\sum_{k=0}^{n-3}\log^k(y)}{y^2}\mathrm{d}y.
\end{split}
\end{equation}
We work with the last two terms in \Cref{eq:int-3-terms} separately. First,
\begin{equation}
\label{eq:int-mid-term-upper-bound}
\begin{split}
\frac{1}{z}\int_{2}^z D_n \frac{\sum_{k=0}^{n-3}\log^k(y)}{y}\mathrm{d}y & = \frac{D_n}{z} \left(\sum_{k=0}^{n-3}\frac{\log^{k+1}(y)}{k+1}\right)\Bigg|_{2}^z \\&\leqslant \frac{D_n}{z} \sum_{k=0}^{n-3}\frac{\log^{k+1}(z)}{k+1} \\&\leqslant D_n \frac{\sum_{k=0}^{n-3}\log^{k+1}(z)}{z}.
\end{split}
\end{equation}
Next, by \cite[Equation~2.722]{MR2360010}, for $n\geqslant0$,
\begin{align}
\int \frac{\log^n(x)}{x^2} \mathrm{d}x  =  -\frac{\sum_{m=0}^n (n-m)!{n \choose k}\log^{m}(x)} {x},
\end{align}
therefore
\begin{equation}
\label{eq:int-long-term-upper-bound}
\begin{split}
\int_z^\infty D_n\frac{\sum_{k=0}^{n-3}\log^k(y)}{y^2}\mathrm{d}y \leqslant  D_n \sum_{k=0}^{n-3} k!\frac{\sum_{m=0}^{k}\log^m(z)}{z} \leqslant D_n(n-2)!\frac{\sum_{k=0}^{n-3}\log^k(z)}{z}.
\end{split}
\end{equation}
Putting \Cref{eq:int-3-terms,eq:int-mid-term-upper-bound,eq:int-long-term-upper-bound} together, if we let
\begin{align*}
D_{n+1} = C_n+D_n+D_n(n-2)!
\end{align*}
and
\begin{align*}
C_{n+1} = \max_{[1,2]}\,\, D_{n+1} \frac{\sum_{k=0}^{n-2}\log^k(y)}{y},
\end{align*}
then 
\begin{equation}
\begin{split}
    \int_1^\infty \frac{V_{\DD_n}(y)}{y+z} \mathrm{d}y  \leqslant C_{n+1},& \quad1\leqslant y\leqslant 2,
    \\
    \int_1^\infty \frac{V_{\DD_n}(y)}{y+z} \mathrm{d}y \leqslant D_{n+1} \frac{\sum_{k=0}^{n-2}\log^k(y)}{y},& \quad 2\leqslant y < \infty.
\end{split}
\end{equation}
Hence, by the Fubini-Tonelli theorem, $\int_1^\infty \frac{V_{\DD_n}(y)}{y+z} \mathrm{d}y = V_{\DD_{n+1}}(z)$, which concludes the proof.
\end{proof}

\begin{proof}[Proof of \Cref{prop:fubini-tonelli}] The first claim follows immediately from \Cref{lem:fubini-tonelli}. Next, 
by \cite[Equation~1.9.7]{MR2254107}, 
\begin{align}
    P_{-1/2+i\xi}(y) = O(y^{-1/2})
\end{align}
as $y \to \infty$. Let $a(\xi)\geqslant 2,C(\xi)>0$ be such that
\begin{equation*}
\begin{split}
|P_{-1/2+i\xi}(y)|\leqslant& C(\xi), \quad\,\,\, 1\leqslant y\leqslant a(\xi), \\
|P_{-1/2+i\xi}(y)|\leqslant& C(\xi)y^{-1/2}, \quad y\geqslant a(\xi).
\end{split}
\end{equation*}
Then by \Cref{lem:fubini-tonelli},
\[
\begin{aligned}
&\int_1^\infty |V_{\DD_n}(y)\cdot P_{-1/2+i\xi}(y)|\,\mathrm{d}y \\
\smash{\mathllap{=}}\, &\int_1^2 |V_{\DD_n}(y)\cdot P_{-1/2+i\xi}(y)|\,\mathrm{d}y+\int_2^{a(\xi)} |V_{\DD_n}(y)\cdot P_{-1/2+i\xi}(y)|\,\mathrm{d}y\\&+\int_{a(\xi)}^{\infty} |V_{\DD_n}(y)\cdot P_{-1/2+i\xi}(y)|\,\mathrm{d}y\\
\smash{\mathllap{\leqslant}}\, &  C_n\cdot C(\xi) + \int_2^{a(\xi)} D_n C(\xi)\frac{\sum_{k=0}^{n-3}\log^k(y)}{y}\mathrm{d}y +\int_{a(\xi)}^\infty D_n C(\xi)\frac{\sum_{k=0}^{n-3}\log^k(y)}{y^{3/2}}\mathrm{d}y \\
\smash{\mathllap{<}}\, & \infty.
\end{aligned}
\]

\end{proof}

\bibliographystyle{plain}
\bibliography{bibliography}

\begin{thebibliography}{10}

\bibitem{zbMATH06877694}
A.~Alexandrov.
\newblock Open intersection numbers, matrix models and {MKP} hierarchy.
\newblock {\em J. High Energy Phys.}, 2015(3):14, 2015.
\newblock Id/No 42.

\bibitem{zbMATH06824078}
Alexander Alexandrov, Alexandr Buryak, and Ran~J. Tessler.
\newblock Refined open intersection numbers and the {Kontsevich}-{Penner} matrix model.
\newblock {\em J. High Energy Phys.}, 2017(3):41, 2017.
\newblock Id/No 123.

\bibitem{zbMATH07023773}
Marco Bertola and Giulio Ruzza.
\newblock The {Kontsevich}-{Penner} matrix integral, isomonodromic tau functions and open intersection numbers.
\newblock {\em Ann. Henri Poincar{\'e}}, 20(2):393--443, 2019.

\bibitem{zbMATH06576752}
A.~Buryak.
\newblock Open intersection numbers and the wave function of the {KdV} hierarchy.
\newblock {\em Mosc. Math. J.}, 16(1):27--44, 2016.

\bibitem{zbMATH06497098}
Alexandr Buryak.
\newblock Equivalence of the open {KdV} and the open {Virasoro} equations for the moduli space of {Riemann} surfaces with boundary.
\newblock {\em Lett. Math. Phys.}, 105(10):1427--1448, 2015.

\bibitem{buryak2017matrix}
Alexandr Buryak and Ran~J Tessler.
\newblock Matrix models and a proof of the open analog of {Witten’s} conjecture.
\newblock {\em Communications in Mathematical Physics}, 353(3):1299--1328, 2017.

\bibitem{chekhov2024}
Leonid~O. Chekhov.
\newblock Fool's crowns, trumpets, and {S}chwarzian.
\newblock {\em arXiv preprint arXiv:2411.03913}, 2024.

\bibitem{MR460128}
Louis Comtet.
\newblock {\em Advanced combinatorics}.
\newblock D. Reidel Publishing Co., Dordrecht, enlarged edition, 1974.
\newblock The art of finite and infinite expansions.

\bibitem{batemanHTF1}
A.~Erd{\'e}lyi, W.~Magnus, F.~Oberhettinger, and F.~G. Tricomi.
\newblock {\em Higher Transcendental Functions. Vol. I}.
\newblock Bateman Manuscript Project. McGraw--Hill, New York, 1953.

\bibitem{galkin_apery}
Sergey Galkin.
\newblock Ap{\'e}ry constants of homogeneous varieties.
\newblock Preprint, {arXiv}:1604.04652 [math.{NT}] (2016), 2016.

\bibitem{MR2254107}
H.-J. Glaeske, A.~P. Prudnikov, and K.~A. Sk\`ornik.
\newblock {\em Operational calculus and related topics}, volume~10 of {\em Analytical Methods and Special Functions}.
\newblock Chapman \& Hall/CRC, Boca Raton, FL, 2006.

\bibitem{goncharovsun}
Alexander~B. Goncharov and Zhe Sun.
\newblock Exponential volumes of moduli spaces of hyperbolic surfaces.
\newblock Preprint, {arXiv}:2411.01615 [math.{AG}] (2024), 2024.

\bibitem{MR2360010}
I.~S. Gradshteyn and I.~M. Ryzhik.
\newblock {\em Table of integrals, series, and products}.
\newblock Elsevier/Academic Press, Amsterdam, seventh edition, 2007.
\newblock Translated from the Russian, Translation edited and with a preface by Alan Jeffrey and Daniel Zwillinger, With one CD-ROM (Windows, Macintosh and UNIX).

\bibitem{huang_thesis}
Yi~Huang.
\newblock {\em Moduli Spaces of Surfaces}.
\newblock PhD thesis, The University of Melbourne, June 2014.

\bibitem{MO402985}
Joe.
\newblock A multiple integral that seems related to the zeta function at even integers.
\newblock \url{https://mathoverflow.net/questions/402985/a-multiple-integral-that-seems-related-to-the-zeta-function-at-even-integers}.
\newblock MathOverflow, accessed September 18, 2025.

\bibitem{MR1171758}
Maxim Kontsevich.
\newblock Intersection theory on the moduli space of curves and the matrix {A}iry function.
\newblock {\em Comm. Math. Phys.}, 147(1):1--23, 1992.

\bibitem{MR350075}
N.~N. Lebedev.
\newblock {\em Special functions and their applications}.
\newblock Dover Publications, Inc., New York, revised edition, 1972.
\newblock Unabridged and corrected republication.

\bibitem{martelli2022introductiongeometrictopology}
Bruno Martelli.
\newblock An introduction to geometric topology.
\newblock {\em arXiv preprint arXiv:1610.02592}, 2016.

\bibitem{MR1510098}
F.~G. Mehler.
\newblock Ueber eine mit den {K}ugel- und {C}ylinderfunctionen verwandte {F}unction und ihre {A}nwendung in der {T}heorie der {E}lektricit\"atsvertheilung.
\newblock {\em Math. Ann.}, 18(2):161--194, 1881.

\bibitem{mirz_simp}
Maryam Mirzakhani.
\newblock Simple geodesics and {W}eil-{P}etersson volumes of moduli spaces of bordered {R}iemann surfaces.
\newblock {\em Invent. Math.}, 167(1):179--222, 2007.

\bibitem{zbMATH07951658}
Rahul Pandharipande, Jake~P. Solomon, and Ran~J. Tessler.
\newblock Intersection theory on moduli of disks, open {KdV} and {Virasoro}.
\newblock {\em Geom. Topol.}, 28(6):2483--2567, 2024.

\bibitem{pennercoords}
R.~C. Penner.
\newblock The decorated {T}eichm\"uller space of punctured surfaces.
\newblock {\em Comm. Math. Phys.}, 113(2):299--339, 1987.

\bibitem{pennerbook}
Robert~C. Penner.
\newblock {\em Decorated {Teichm{\"u}ller} theory}.
\newblock QGM Master Cl. Ser. Z{\"u}rich: European Mathematical Society (EMS), 2012.

\bibitem{MR2868112}
Richard~P. Stanley.
\newblock {\em Enumerative combinatorics. {V}olume 1}, volume~49 of {\em Cambridge Studies in Advanced Mathematics}.
\newblock Cambridge University Press, Cambridge, second edition, 2012.

\bibitem{weil}
Andr\'{e} Weil.
\newblock Modules des surfaces de {R}iemann.
\newblock {\em S\'{e}minaire Bourbaki}, 4:413--419, 1956-1958.

\bibitem{wikipedia:dirichlet_beta}
{Wikipedia contributors}.
\newblock Dirichlet beta function --- {W}ikipedia{,} the free encyclopedia.
\newblock \url{https://en.wikipedia.org/wiki/Dirichlet_beta_function}, 2025.
\newblock [Online; accessed 9-September-2025].

\bibitem{witten_conjecture}
Edward Witten.
\newblock Two-dimensional gravity and intersection theory on moduli space.
\newblock In {\em Surveys in differential geometry ({C}ambridge, {MA}, 1990)}, pages 243--310. Lehigh Univ., Bethlehem, PA, 1991.

\bibitem{wolpertwpform}
Scott Wolpert.
\newblock On the {W}eil-{P}etersson geometry of the moduli space of curves.
\newblock {\em Amer. J. Math.}, 107(4):969--997, 1985.

\end{thebibliography}
\end{document}